\documentclass{article}
\usepackage{amsmath}
\usepackage{amssymb}
\usepackage{mathtools}
\usepackage{enumerate}
\usepackage{graphicx}
\usepackage{float}
\usepackage{overpic}
\usepackage{subfigure}
\usepackage{caption}
\usepackage{relsize}
\usepackage{geometry}
\usepackage{algorithm}
\usepackage{algorithmic}
\usepackage{multirow}
\usepackage{stmaryrd}
\usepackage{tikz}
\usetikzlibrary{decorations.pathmorphing}
\usepackage{epstopdf}
\newcommand\dd{\,\mathrm{d}}
\newcommand\ii{\,\mathrm{i}}

\newcommand\td{\tilde}

\begin{document}

\title{Additive Sweeping Preconditioner for the Helmholtz Equation}
\author{Fei Liu$^\sharp$ and Lexing Ying$^{\dagger\sharp}$\\
  $\dagger$ Department of Mathematics, Stanford University\\
  $\sharp$ Institute for Computational and Mathematical Engineering, Stanford University
}
\date{Mar. 2016}
\maketitle

\begin{abstract}
  We introduce a new additive sweeping preconditioner for the
  Helmholtz equation based on the perfect matched layer (PML). This
  method divides the domain of interest into thin layers and proposes
  a new transmission condition between the subdomains where the
  emphasis is on the boundary values of the intermediate waves. This
  approach can be viewed as an effective approximation of an additive
  decomposition of the solution operator. When combined with the
  standard GMRES solver, the iteration number is essentially
  independent of the frequency. Several numerical examples are tested
  to show the efficiency of this new approach.
\end{abstract}


\section{Introduction}
Let the domain of interest be $D=(0,1)^d$ where $d=2,3$. The Helmholtz
equation is
\begin{equation}
  \Delta u(x)+\dfrac{\omega^2}{c^2(x)}u(x)=f(x),\quad \forall x\in D,
  \label{eq:helm}
\end{equation}
where $u(x)$ is the time-independent wave field generated by the
time-independent force $f(x)$, $\omega$ is the angular frequency and
$c(x)$ is the velocity field. Commonly used boundary conditions are
the approximations of the Sommerfeld radiation condition. By rescaling the system, we assume $c_{\min}\le c(x)\le c_{\max}$ where $c_{\min}$
and $c_{\max}$ are of $\Theta(1)$. Then $\omega/(2\pi)$ is the typical
wave number and $\lambda=2\pi/\omega$ is the typical wavelength.

Solving the equation numerically is challenging in high frequency
settings for two reasons. First, in most applications, the equation is
discretized with at least a constant number of points per wavelength,
which makes the number of points in each direction $n=\Omega(\omega)$
and the total degree of freedom $N=n^d=\Omega(\omega^d)$ very
large. Second, the system is highly indefinite and has a very
oscillatory Green's function, which makes most of the classical
iterative methods no longer effective.

There has been a sequence of papers on developing iterative methods
for solving \eqref{eq:helm}. The AILU method by Gander and Nataf
\cite{ailu} is the first to use the incomplete LU factorization to
precondition the equation. Engquist and Ying \cite{sweephmf,sweeppml}
developed a series of sweeping preconditioners based on approximating
the inverse of the Schur complements in the LDU factorization and
obtained essentially $\omega$-independent iteration numbers. In
\cite{stolk2013domaindecomp}, Stolk proposed a domain decomposition
method based on the PML which constructs delicate transmission
conditions between the subdomains by considering the ``pulses''
generated by the intermediate waves. In \cite{vion2014doublesweep},
Vion and Geuzaine proposed a double sweep preconditioner based on the
Dirichlet-to-Neumann (DtN) map and several numerical simulations of
the DtN map were compared. In
\cite{chen2013sourcetrans,chen2013sourcetrans2}, Chen and Xiang
introduced a source transfer domain decomposition method which
emphasizes on transferring the sources between the subdomains. In
\cite{demanet}, Zepeda-N{\'u}{\~n}ez and Demanet developed a novel
domain decomposition method for the 2D case by pairing up the waves
and their normal derivatives at the boundary of the subdomains and
splitting the transmission of the waves into two directions. Most
recently in \cite{Liu2015}, Liu and Ying proposed a recursive sweeping
preconditioner for 3D Helmholtz problems. Other progresses includes
\cite{parallelsweep,sweepem,sweepemfem,sweepspectral} and we refer to
\cite{advances} by Erlangga and \cite{why} by Ernst and Gander for a
complete discussion.

Inspired by \cite{stolk2013domaindecomp} and these previous
approaches, we propose a new domain decomposition method in this paper
which shares some similarities with
\cite{sweeppml,stolk2013domaindecomp}. The novelty of this new
approach is that the transmission conditions are built with the
boundary values of the intermediate waves directly. For each wave
field on the subdomains, we divide it into three parts -- the waves
generated by the force to the left of the subdomain, to the right of
the subdomain, and within the subdomain itself. This corresponds to an
$L+D+U$ decomposition of the Green's matrix $G$ as the sum of its
lower triangular part, upper triangular part and diagonal part. This
is why we call this new preconditioner the additive sweeping
preconditioner.

The rest of this paper is organized as follows. First in Section
\ref{sec:1D} we use the 1D case to illustrate the idea of the
method. Then in Section \ref{sec:2D} we introduce the preconditioner
in 2D and present the 2D numerical results. Section
\ref{sec:3D} discusses the 3D case.  Conclusions and some future
directions are provided in Section \ref{sec:Conclusion}.

\section{1D Illustration}
\label{sec:1D}

We use the PML\cite{berenger1994pml,chew1994pml,johnson2008pmlnotes}
to simulate the Sommerfeld condition. The PML introduces the auxiliary
functions
\begin{align*}
\sigma(x) &:=
\begin{dcases}
\dfrac{C}{\eta}\left(\dfrac{x-\eta}{\eta}\right)^2,\quad &x\in[0,\eta),\\
0,\quad &x\in[\eta,1-\eta],\\
\dfrac{C}{\eta}\left(\dfrac{x-1+\eta}{\eta}\right)^2,\quad &x\in(1-\eta,1],
\end{dcases}\\
s(x) &:= \left(1+\ii\dfrac{\sigma(x)}{\omega}\right)^{-1},
\end{align*}
where C is an appropriate positive constant independent of $\omega$,
and $\eta$ is the PML width which is typically around one wavelength.

The Helmholtz equation with PML in 1D is
\begin{equation*}
\begin{dcases}
\left((s(x)\dfrac{\dd}{\dd x})^2+\dfrac{\omega^2}{c^2(x)}\right)u(x)=f(x),\quad \forall x\in (0,1),\\
u(0)=0,\\
u(1)=0. 
\end{dcases}
\end{equation*}
We discretize the system with step size $h=1/(n+1)$, then $n$
is the degree of freedom. With the standard central difference
numerical scheme the discretized equation is
\begin{equation}
\label{eqn:1D}
\dfrac{s_{i}}{h}\left(\dfrac{s_{i+1/2}}{h}(u_{i+1}-u_{i})-\dfrac{s_{i-1/2}}{h}(u_{i}-u_{i-1})\right)+\dfrac{\omega^2}{c_i^2}u_{i}=f_{i}, \quad \forall 1\le i\le n,
\end{equation}
where the subscript $i$ means that the corresponding function is evaluated at $x=ih$.

We denote Equation \eqref{eqn:1D} as $A\pmb u=\pmb f$, where $\pmb u$ and $\pmb f$ are the discrete array of the wave field and the
force
\begin{align*}
\pmb u:=[u_1,\dots,u_n]^T,\quad \pmb f:=[f_1,\dots,f_n]^T. 
\end{align*}
In 1D, $A$ is tridiagonal and Equation \eqref{eqn:1D} can be solved
without any difficulty. However, here we are aiming at an approach
which can be generalized to higher dimensions so the rest of this
section takes another point of view to solve \eqref{eqn:1D} instead of
exploiting the sparsity structure of $A$ directly.

With the Green's matrix $G=A^{-1}$, $\pmb u$ can be written
as $\pmb u=G\pmb f$. Now let us divide the discrete grid into $m$ parts. We
assume that $\eta=\gamma h$ and $n=2\gamma+mb-2$ where $\gamma$ and $b$ are some
small constants and $m$ is comparable to $n$, and we define
\begin{align*}
X_1 &:= \{ih:1 \le i \le \gamma + b-1 \}, \\
X_p &:= \{ih:\gamma + (p-1) b \le i \le \gamma + pb-1 \},\quad p=2,\dots, m-1, \\
X_m &:= \{ih:\gamma + (m-1) b \le i \le 2 \gamma + mb-2 \}, 
\end{align*}
which means, $X_1$ is the leftmost part containing the left PML of the
original problem and a small piece of grid with $b$ points, $X_m$ is
the rightmost part containing the right PML and a grid of $b$ points,
and $X_p, p=2,\dots,m-1$ are the middle parts each of which contains
$b$ points. $\pmb u_p$ and $\pmb f_p$ are defined as the restrictions
of $\pmb u$ and $\pmb f$ on $X_p$ for $p=1,\dots,m$, respectively,
\begin{align*}
\pmb u_1 &:=  [u_1,\dots,u_{\gamma+b-1}]^T,\\
\pmb u_p &:=  [u_{\gamma+(p-1)b},\dots,u_{\gamma + pb-1}]^T,\quad p=2,\dots,m-1,\\
\pmb u_m &:=  [u_{\gamma + (m-1) b},\dots,u_{2 \gamma + m b -2}]^T,\\
\pmb f_1 &:=  [f_1,\dots,f_{\gamma+b-1}]^T,\\
\pmb f_p &:=  [f_{\gamma+(p-1)b},\dots,f_{\gamma + pb-1}]^T,\quad p=2,\dots,m-1,\\
\pmb f_m &:=  [f_{\gamma + (m-1) b},\dots,f_{2 \gamma + m b -2}]^T.
\end{align*}
Then $u=Gf$ can be written as
\begin{align*}
\begin{bmatrix}
\pmb u_1\\\pmb u_2\\ \vdots\\\pmb u_m
\end{bmatrix}
=
\begin{bmatrix}
G_{1,1}&G_{1,2}&\ldots&G_{1,m}\\
G_{2,1}&G_{2,2}&\ldots&G_{2,m}\\
\vdots&\vdots&&\vdots\\
G_{m,1}&G_{m,2}&\ldots&G_{m,m}
\end{bmatrix}
\begin{bmatrix}
\pmb f_1\\\pmb f_2\\ \vdots\\\pmb f_m
\end{bmatrix}.
\end{align*}

By introducing $\pmb u_{p,q}:=G_{p,q}\pmb f_q$ for $1\le p,q\le m$,
one can write $\pmb u_p=\sum_{q=1}^m \pmb u_{p,q}$. The physical
meaning of $\pmb u_{p,q}$ is the contribution of the force $\pmb f_q$
defined on the grid $X_q$ acting upon the grid $X_p$. If we know the
matrix $G$, the computation of $\pmb u_{p,q}$ can be carried out
directly. However, computing $G$, or even applying $G$ to the vector
$\pmb f$, is computationally expensive. The additive sweeping method
circumvent this difficulty by approximating the blocks of $G$
sequentially and the idea works in higher dimensions. In what follows,
we shall use $\td{\pmb u}_{p,q}$ to denote the approximations of $\pmb
u_{p,q}$.

\subsection{Approximating $\pmb u_{p,q}$ with auxiliary PMLs}

\subsubsection{Wave generated by $\pmb f_1$}

\begin{figure}[h!]
  \centering
    
    \begin{tikzpicture}
      [x=2cm,y=2cm,>=latex]
      \draw
      (0.0 , 0.0)--(7.0 , 0.0)
      
      (0.0 , 1.5)--(2.0 , 1.5)
      (2.5 , 1.5)--(4.5 , 1.5)
      (5.0 , 1.5)--(7.0 , 1.5)
      
      (0.0 , 1.0)--(1.5 , 1.0)
      (2.5 , 1.0)--(4.0 , 1.0)
            
      (3.0 , 2.0)--(4.5 , 2.0)
      (5.5 , 2.0)--(7.0 , 2.0)
      ;
      
      \draw
      (0.0 , -0.1)--(0.0 , 0.1)
      (1.5 , -0.1)--(1.5 , 0.1)
      (3.0 , -0.1)--(3.0 , 0.1)
      (4.0 , -0.1)--(4.0 , 0.1)
      (5.5 , -0.1)--(5.5 , 0.1)
      (7.0 , -0.1)--(7.0 , 0.1)
      
      (0.0 , 1.0 - 0.1)--(0.0 , 1.0 + 0.1)
      (1.5 , 1.0 - 0.1)--(1.5 , 1.0 + 0.1)
      (2.5 , 1.0 - 0.1)--(2.5 , 1.0 + 0.1)
      (4.0 , 1.0 - 0.1)--(4.0 , 1.0 + 0.1)
      
      (3.0 , 2.0 - 0.1)--(3.0 , 2.0 + 0.1)
      (4.5 , 2.0 - 0.1)--(4.5 , 2.0 + 0.1)
      (5.5 , 2.0 - 0.1)--(5.5 , 2.0 + 0.1)
      (7.0 , 2.0 - 0.1)--(7.0 , 2.0 + 0.1)
      
      (0.0 , 1.5 - 0.1)--(0.0 , 1.5 + 0.1)
      (2.0 , 1.5 - 0.1)--(2.0 , 1.5 + 0.1)
      (2.5 , 1.5 - 0.1)--(2.5 , 1.5 + 0.1)
      (4.5 , 1.5 - 0.1)--(4.5 , 1.5 + 0.1)
      (5.0 , 1.5 - 0.1)--(5.0 , 1.5 + 0.1)
      (7.0 , 1.5 - 0.1)--(7.0 , 1.5 + 0.1)
      ;
      
      \draw
      (0.0 , 0.0)node[anchor=south west]{PML}
      (0.0 , 1.0)node[anchor=south west]{PML}
      (0.0 , 1.5)node[anchor=south west]{PML}
      (2.5 , 1.0)node[anchor=south west]{PML}
      (2.5 , 1.5)node[anchor=south west]{PML}
      (5.0 , 1.5)node[anchor=south west]{PML}
      (7.0 , 0.0)node[anchor=south east]{PML}
      (7.0 , 1.5)node[anchor=south east]{PML}
      (7.0 , 2.0)node[anchor=south east]{PML}
      (4.5 , 1.5)node[anchor=south east]{PML}
      (4.5 , 2.0)node[anchor=south east]{PML}
      (2.0 , 1.5)node[anchor=south east]{PML}
      ;
      
      \draw
      (1.0 , 0.0)node[anchor = north]{$X_1$}
      (1.0 , 1.0)node[anchor = north]{$X_1^L$}
      (1.0 , 1.5)node[anchor = north]{$X_1^M$}
      
      (6.0 , 0.0)node[anchor = north]{$X_m$}
      (6.0 , 2.0)node[anchor = north]{$X_m^R$}
      (6.0 , 1.5)node[anchor = north]{$X_m^M$}
      
      (3.5 , 0.0)node[anchor = north]{$X_p$}
      (3.5 , 1.0)node[anchor = north]{$X_p^L$}
      (3.5 , 1.5)node[anchor = north]{$X_p^M$}
      (3.5 , 2.0)node[anchor = north]{$X_p^R$}
      ;
      
      \draw
      (2.25 , 0.0)node[fill=white]{$\ldots$}
      (4.75 , 0.0)node[fill=white]{$\ldots$}
      ;
    \end{tikzpicture}
  
  \caption{This figure shows how the grids $X_p$ are extended with auxiliary PMLs.}
  \label{fig:Xp}
\end{figure}
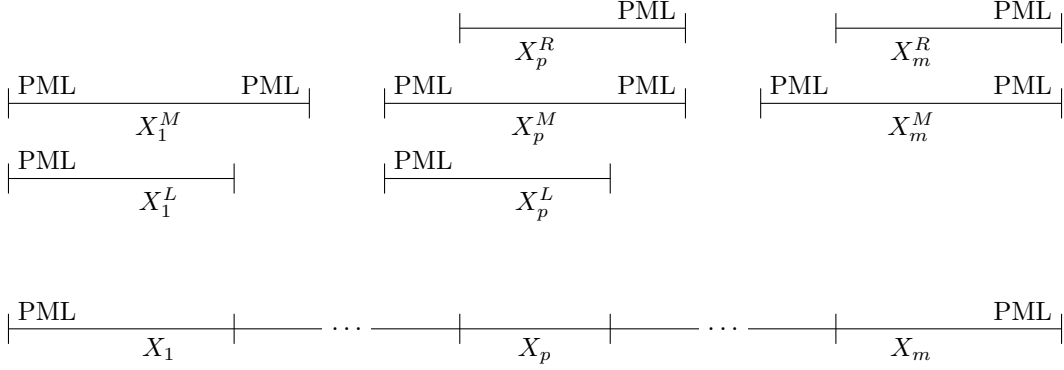

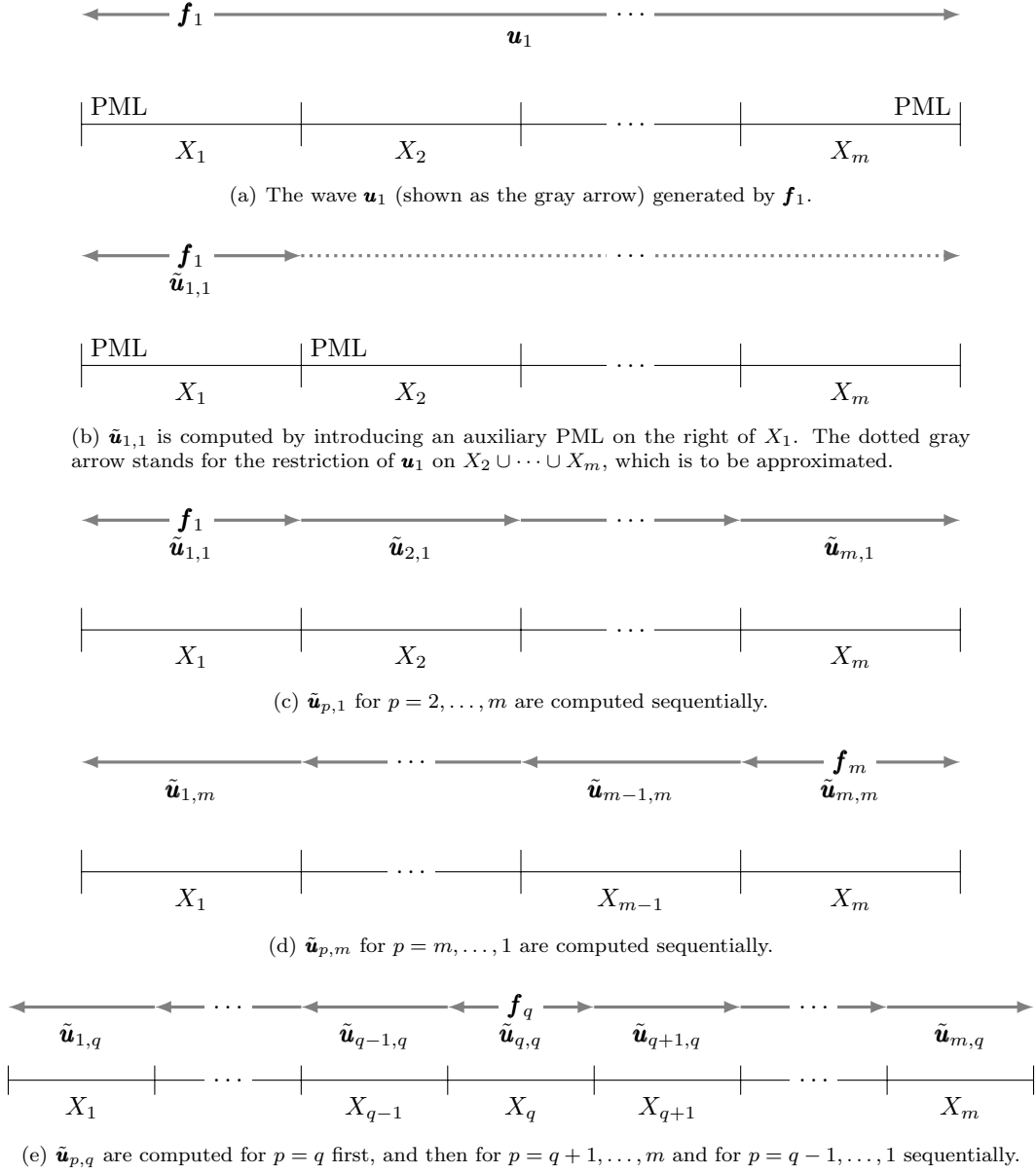
\begin{figure}[h!]
  \centering
  \subfigure[The wave $\pmb u_1$ (shown as the gray arrow) generated by $\pmb f_1$.]{
    \begin{tikzpicture}
      [x=3cm,y=3cm,>=latex]
      \draw[step=1](0,-0.1)grid(4,0.1);
      \draw(0.5,0)node[below=3]{$X_1$}(1.5,0)node[below=3]{$X_2$}(3.5,0)node[below=3]{$X_m$};
      \draw[<->,gray,very thick](0,0.5)--(4,0.5);
      \draw(0.5,0.5)node[fill=white]{$\pmb f_1$};
      \draw(2.0,0.5)node[below=3]{$\pmb u_1$};
      \draw(2.5,0.5)node[fill=white]{$\ldots$}(2.5,0)node[fill=white]{$\ldots$};
      \draw(0,0)node[anchor=south west]{PML}(4,0)node[anchor=south east]{PML};
    \end{tikzpicture}
  }
  \subfigure[$\td{\pmb u}_{1,1}$ is computed by introducing an auxiliary PML on the right of $X_1$. The dotted gray arrow stands for the restriction of ${\pmb u}_1$ on $X_2\cup\dots \cup X_m$, which is to be approximated.]{
    \label{fig:approxU11}
    \begin{tikzpicture}
      [x=3cm,y=3cm,>=latex]
      \draw[step=1](0,-0.1)grid(4,0.1);
      \draw(0.5,0)node[below=3]{$X_1$}(1.5,0)node[below=3]{$X_2$}(3.5,0)node[below=3]{$X_m$};
      \draw[<->,gray,very thick](0,0.5)--(1,0.5);
      \draw[->,dotted,gray,very thick](1,0.5)--(4,0.5);
      \draw(0.5,0.5)node[fill=white]{$\pmb f_1$};
      \draw(0.5,0.5)node[below=3]{$\td{\pmb u}_{1,1}$};
      \draw(2.5,0.5)node[fill=white]{$\ldots$}(2.5,0)node[fill=white]{$\ldots$};
      \draw(0,0)node[anchor=south west]{PML}(1,0)node[anchor=south west]{PML};
    \end{tikzpicture}
  }
  \subfigure[$\td {\pmb u}_{p,1}$ for $p=2,\dots,m$ are computed sequentially.]{
    \label{fig:approxUp1}
    \begin{tikzpicture}
      [x=3cm,y=3cm,>=latex]
      \draw[step=1](0,-0.1)grid(4,0.1);
      \draw(0.5,0)node[below=3]{$X_1$}(1.5,0)node[below=3]{$X_2$}(3.5,0)node[below=3]{$X_m$};
      \draw[<->,gray,very thick](0,0.5)--(1,0.5);
      \foreach \x in{1,2,3}
               {
                 \draw[->,gray,very thick](\x,0.5)--(\x+1,0.5);
               }
               \draw(0.5,0.5)node[fill=white]{$\pmb f_1$};
               \draw(0.5,0.5)node[below=3]{$\td{\pmb u}_{1,1}$}(1.5,0.5)node[below=3]{$\td {\pmb u}_{2,1}$}(3.5,0.5)node[below=3]{$\td {\pmb u}_{m,1}$};
               \draw(2.5,0.5)node[fill=white]{$\ldots$}(2.5,0)node[fill=white]{$\ldots$};
    \end{tikzpicture}
  }
  \subfigure[$\td {\pmb u}_{p,m}$ for $p=m,\dots,1$ are computed sequentially.]
  {
  \label{fig:approxUpm}
  \begin{tikzpicture}
    [x=3cm,y=3cm,>=latex]
    \draw[step=1](0,-0.1)grid(4,0.1);
    \draw(0.5,0)node[below=3]{$X_1$}(2.5,0)node[below=3]{$X_{m-1}$}(3.5,0)node[below=3]{$X_m$};
    \draw[<->,gray,very thick](3,0.5)--(4,0.5);
    \foreach \x in{3,2,1}
             {
               \draw[->,gray,very thick](\x,0.5)--(\x-1,0.5);
             }
             \draw(3.5,0.5)node[fill=white]{$\pmb f_m$};
             \draw(0.5,0.5)node[below=3]{$\td {\pmb u}_{1,m}$}(2.5,0.5)node[below=3]{$\td {\pmb u}_{m-1,m}$}(3.5,0.5)node[below=3]{$\td {\pmb u}_{m,m}$};
             \draw(1.5,0.5)node[fill=white]{$\ldots$}(1.5,0)node[fill=white]{$\ldots$};
  \end{tikzpicture}
  }
  \subfigure[$\td {\pmb u}_{p,q}$ are computed for $p=q$ first, and then for $p=q+1,\dots,m$ and for $p=q-1,\dots,1$ sequentially.]
  {
  \label{fig:approxUpq}
  \begin{tikzpicture}
    [x=2cm,y=2cm,>=latex]
    \draw[step=1](0,-0.1)grid(7,0.1);
    \draw(0.5,0)node[below=3]{$X_1$}(2.5,0)node[below=3]{$X_{q-1}$}(3.5,0)node[below=3]{$X_q$}(4.5,0)node[below=3]{$X_{q+1}$}(6.5,0)node[below=3]{$X_m$};
    \draw[<->,gray,very thick](3,0.5)--(4,0.5);
    \foreach \x in{3,2,1}
             {
               \draw[->,gray,very thick](\x,0.5)--(\x-1,0.5);
               \draw[->,gray,very thick](\x+3,0.5)--(\x+4,0.5);
             }
             \draw(3.5,0.5)node[fill=white]{$\pmb f_q$};
             \draw(0.5,0.5)node[below=3]{$\td {\pmb u}_{1,q}$}(2.5,0.5)node[below=3]{$\td {\pmb u}_{q-1,q}$}(3.5,0.5)node[below=3]{$\td {\pmb u}_{q,q}$}(4.5,0.5)node[below=3]{$\td {\pmb u}_{q+1,q}$}(6.5,0.5)node[below=3]{$\td {\pmb u}_{m,q}$};
             \draw(1.5,0.5)node[fill=white]{$\ldots$}(1.5,0)node[fill=white]{$\ldots$}(5.5,0.5)node[fill=white]{$\ldots$}(5.5,0)node[fill=white]{$\ldots$};
  \end{tikzpicture}
  }
  \caption{This figure shows how $\td {\pmb u}_{p,q}$ are generated. The direction of the arrows indicates the computing orders of the approximating waves.}
\end{figure}

The components ${\pmb u}_{p,1}$ for $p=1,\dots,m$ can be regarded as a
sequence of right-going waves generated by $\pmb f_1$. Note that the
boundary condition of the system is the approximated Sommerfeld
condition. If we assume that the reflection during the transmission of
the wave is negligible, then, to approximate $\pmb u_{1,1}$, we can
simply put an artificial PML on the right of the grid $X_1$ to solve a
much smaller problem, since the domain of interest here is only $X_1$
(see Figure \ref{fig:approxU11}). To be precise, we define
\begin{align*}
\sigma_1^{M}(x) &:=  
\begin{dcases}
\dfrac{C}{\eta}\left(\dfrac{x-\eta}{\eta}\right)^2,\quad &x\in[0,\eta),\\
0,\quad &x\in[\eta,\eta+(b-1)h],\\
\dfrac{C}{\eta}\left(\dfrac{x-(\eta+(b-1)h)}{\eta}\right)^2,\quad &x\in(\eta+(b-1)h,2\eta+(b-1)h],
\end{dcases}\\
s_1^{M}(x) &:=  \left(1+\ii\dfrac{\sigma_1^{M}(x)}{\omega}\right)^{-1}. 
\end{align*}
We consider a subproblem on the auxiliary domain $D_1^M := (0,2\eta + (b - 1)h)$
\begin{align*}
\begin{dcases}
\left((s_1^{M}(x)\dfrac{\dd}{\dd x})^2+\dfrac{\omega^2}{c^2(x)}\right)v(x)=g(x),\quad &\forall x\in D_1^M,\\
v(x)=0, \quad & \forall x\in \partial D_1^M.
\end{dcases}
\end{align*}
With the same discrete numerical scheme and step size $h$, we have the
corresponding discrete system $H_1^{M} \pmb v=\pmb g$ on the extended grid
\begin{align*}
X_1^{M}:=\{ih :1\le i\le 2\gamma+b-2\}. 
\end{align*}
Figure \ref{fig:Xp} shows a graphical view of $X_1^M$, as well as other extended grids which we will see later.

With the discrete system $H_1^M \pmb v = \pmb g$, we can define an operator $\td{G}_{1}^{M}:\pmb y\to \pmb z$, which is an approximation of
$G_{1,1}$, by the following:
\begin{enumerate}
\item 
  Introduce a vector $\pmb g$ defined on $X_1^{M}$ by setting
  $\pmb y$ to $X_1$ and zero everywhere else.
\item 
  Solve $H_1^{M} \pmb v=\pmb g$ on $X_1^{M}$.
\item 
  Set $\pmb z$ as the restriction of $\pmb v$ on $X_1$.
\end{enumerate}
Then $\td{\pmb u}_{1,1}$ can be set as
\begin{align*}
\td{\pmb u}_{1,1}:=\td{G}_{1}^{M}\pmb f_1.
\end{align*}

Once we have computed $\td{\pmb u}_{1,1}$, we can use the right
boundary value of $\td{\pmb u}_{1,1}$ to compute $\td{\pmb u}_{2,1}$
by introducing an auxiliary PML on the right of $X_2$ and solving the
boundary value problem with the left boundary value at
$x=(\gamma+b-1)h$ equal to the right boundary value of $\td{\pmb
  u}_{1,1}$. The same process can be repeated to compute $\td{\pmb
  u}_{p+1,1}$ by exploiting the right boundary value of $\td{\pmb
  u}_{p,1}$ recursively for $p=2,\dots,m-1$ (see Figure
\ref{fig:approxUp1}). In the following context of this section, we
introduce notations $g^{L}, g^{R}$ for a vector array $\pmb
g=[g_1,\dots,g_s]^T$ by
\begin{align*}
g^{L}:=g_1,\quad g^{R}:=g_s,
\end{align*}
where $g^{L}$ and $g^{R}$ should be interpreted as the
leftmost and the rightmost element of the array $\pmb g$.

To formalize the definition of $\td{\pmb u}_{p,1}$ for each
$p=2,\dots,m$, we introduce the auxiliary domain $D_p^R$, which will be defined below, to simulate the right-transmission of the waves. The superscript $R$ means that the auxiliary domain is intended for approximating the right-going waves. The left boundary of $D_p^R$ will be denoted as $\partial^L D_p^R$, on which the boundary value will be used to approximate the wave transmission as we shall see. We also extend $X_p$ with an auxiliary PML on the right to form an extended grid $X_p^R$ (see Figure \ref{fig:Xp}), which corresponds the discretization of $D_p^R$. To be specific, we define
\begin{align*}
D_p^R &:=   (\eta+((p-1)b-1)h,2\eta+(pb-1)h), \\
\partial^L D_p^R  &:=  \{\eta+((p-1)b-1)h\}, \\
  X_p^{R} &:=  \{ih:\gamma+(p-1)b\le i\le 2\gamma+pb-2\}.
\end{align*}
Note that the grid $X_m^{R}$ is $X_m$ itself since
$X_m$ already contains the original right PML region. The purpose to use the notation $X_m^{R}$ is to simplify the description of the algorithm.

For the PML on $D_p^{R}$, we define
\begin{align*}
  \sigma_p^{R}(x) &:=  
  \begin{dcases}
    0,\quad &x\in[\eta+((p-1)b-1)h,\eta+(pb-1)h],\\
    \dfrac{C}{\eta}\left(\dfrac{x-(\eta+(pb-1)h)}{\eta}\right)^2,\quad &x\in(\eta+(pb-1)h,2\eta+(pb-1)h],
  \end{dcases}\\
  s_p^{R} &:=  \left(1+\ii\dfrac{\sigma_p^{R}(x)}{\omega}\right)^{-1}.
\end{align*} 
We consider the following subproblem
\begin{align*}
  \begin{dcases}
    \left((s_p^{R}(x)\dfrac{\dd}{\dd x})^2+\dfrac{\omega^2}{c^2(x)}\right)v(x)=0,\quad &\forall x\in D_p^R,\\
    v(x)=w,\quad & \forall x \in \partial^L D_p^R,\\
    v(x)=0, \quad & \forall x \in \partial D_p^R \setminus\partial^L D_p^R,
  \end{dcases}
\end{align*}
where $w$ is the left boundary value of the unknown $v(x)$. We define $H_p^{R} \pmb v=\pmb g$ as the discretization of this
problem on $X_p^{R}$ where the right-hand side $\pmb g$ is given by
$\pmb g:= (-1/h^2)[w,0,\dots,0]^T$ as a result of the central
discretization. The subproblem $H_p^{R} \pmb v=\pmb g$ for each
$p=2,\dots,m$ induces the approximation operator $\td{G}_{p}^{R}:w \to
\pmb z$ by the following procedure:
\begin{enumerate}
\item
  Set $\pmb g=(-1/h^2)[w,0,\dots,0]^T$. 
\item
  Solve $H_p^{R} \pmb v =\pmb g$ on $X_p^{R}$. 
\item
  Set $\pmb z$ as the restriction of $\pmb v$ on $X_p$.
\end{enumerate}
Then $\td{\pmb u}_{p,1}$ can be defined recursively for $p=2,\dots,m$ by
\begin{align*}
  \td{\pmb u}_{p,1}:=\td{G}_{p}^{R}\td{\pmb u}_{p-1,1}^{R}.
\end{align*}
Note that, the operator $\td{G}_p^{R}$ is not an approximation of the
matrix block $G_{p,1}$, since $\td{G}_p^R$ maps the right boundary
value of $\td{\pmb u}_{p-1,1}$ to $\td{\pmb u}_{p,1}$ while $G_{p,1}$
maps $\pmb f_1$ to $\pmb u_{p,1}$.

\subsubsection{Wave generated by $\pmb f_m$}

The components ${\pmb u}_{p,m}$ for $p=1,\dots,m$ can be regarded as a
sequence of left-going waves generated by $\pmb f_m$. The method for
approximating them is similar to what was done for $\pmb f_1$ (see
Figure \ref{fig:approxUpm}). More specifically, for $\td{\pmb
  u}_{m,m}$ we define
\begin{align*}
  D_m^M  &:=   (1-2\eta-(b-1)h,1),\\
  X_m^{M} &:=  \{ih:(m-1)b+1\le i \le 2\gamma+mb-2\},\\
  \sigma_m^{M}(x) &:=  
  \begin{dcases}
    \dfrac{C}{\eta}\left(\dfrac{x-(1-\eta-(b-1)h)}{\eta}\right)^2,\quad &x\in[1-2\eta-(b-1)h,1-\eta-(b-1)h),\\
      0,\quad &x\in[1-\eta-(b-1)h,1-\eta],\\
      \dfrac{C}{\eta}\left(\dfrac{x-(1-\eta)}{\eta}\right)^2,\quad &x\in(1-\eta,1],
  \end{dcases}\\
  s_m^{M}(x) &:=  \left(1+\ii\dfrac{\sigma_m^{M}(x)}{\omega}\right)^{-1}. 
\end{align*}
We consider the continuous problem
\begin{align*}
  \begin{dcases}
    \left((s_m^{M}(x)\dfrac{\dd}{\dd x})^2+\dfrac{\omega^2}{c^2(x)}\right)v(x)=g(x),\quad & \forall x\in D_m^M,\\
    v(x)=0,\quad & \forall x\in \partial D_m^M,
  \end{dcases}
\end{align*}
and define $H_m^{M}\pmb v=\pmb g$ as its discretization on $X_m^{M}$.
The operator $\td{G}_{m}^{M}: \pmb y\to \pmb z$ can be defined as:
\begin{enumerate}
\item
  Introduce a vector $\pmb g$ defined on $X_m^{M}$ by setting $\pmb y$ to
  $X_m$ and zero everywhere else.
\item
  Solve $H_m^{M} \pmb v=\pmb g$ on $X_m^{M}$. 
\item
  Set $\pmb z$ as the restriction of $\pmb v$ on $X_m$. 
\end{enumerate}
Then
\begin{align*}
  \td{\pmb u}_{m,m}:=\td{G}_{m}^{M}\pmb f_m. 
\end{align*}

For each $\td{\pmb u}_{p,m},p=1,\dots,m-1$, we introduce the auxiliary domain $D_p^L$, the right boundary $\partial ^R D_p^L$, the extended grid $X_p^L$, and the corresponding PML functions $\sigma_p^L(x), s_p^L(x)$ as follows
\begin{align*}
  D_p^L  &:=   ((p-1)bh,\eta+pbh),\\
  \partial ^R D_p^L  &:=  \{\eta + pbh\},\\
  X_p^{L} &:=  \{x_i:(p-1)b+1\le i\le \gamma+pb-1\},\\
  \sigma_p^{L}(x) &:=  
  \begin{dcases}
    \dfrac{C}{\eta}\left(\dfrac{x-(\eta+(p-1)bh)}{\eta}\right)^2,\quad &x\in[(p-1)bh,\eta+(p-1)bh),\\
      0,\quad &x\in[\eta+(p-1)bh,\eta+pbh],
  \end{dcases}\\
  s_p^{L}(x) &:=  \left(1+\ii\dfrac{\sigma_p^{L}(x)}{\omega}\right)^{-1}, 
\end{align*}
and we consider the continuous problem
\begin{align*}
  \begin{dcases}
    \left((s_p^{L}(x)\dfrac{\dd}{\dd
      x})^2+\dfrac{\omega^2}{c^2(x)}\right)v(x)=0,\quad & \forall x\in D_p^L,\\
      v(x) = w,\quad & \forall x \in \partial^R D_p^L,\\
      v(x) = 0,\quad & \forall x \in \partial D_p^L \setminus \partial^R D_p^L,
  \end{dcases}
\end{align*}
where $y$ is the right boundary value of $v(x)$. Let $H_p^{L}\pmb v=\pmb g$ be its discretization on $X_p^{L}$ with
$\pmb g:=(-1/h^2)[0,\dots,0,w]^T$. We introduce the operator
$\td{G}_{p}^{L}:w\mapsto \pmb z$ by:
\begin{enumerate}
\item
  Set $\pmb g=(-1/h^2)[0,\dots,0,w]^T$. 
\item
  Solve $H_p^{L} \pmb v =\pmb g$ on $X_p^{L}$. 
\item
  Set $\pmb z$ as the restriction of $\pmb v$ on $X_p$. 
\end{enumerate}
Then $\td{\pmb u}_{p,m}$ can be defined recursively for
$p=m-1,\dots,1$ by
\begin{align*}
\td{\pmb u}_{p,m}:=\td{G}_{p}^{L}\td{\pmb u}_{p+1,m}^{L}. 
\end{align*}

\subsubsection{Wave generated by $\pmb f_q$ for $q=2,\dots,m-1$}

For each $q$, the components ${\pmb u}_{p,q}$ for $p=1,\dots,m$ can
be regarded as a sequence of left- and right-going waves generated by
$\pmb f_q$ (see Figure \ref{fig:approxUpq}). For $\td{\pmb u}_{q,q}$,
we introduce
\begin{align*}
D_q^M &:=  ((q-1)bh,2\eta+(qb-1)h),\\
X_q^{M} &:=  \{x_i:(q-1)b+1\le i \le 2\gamma+qb-2\},\\
\sigma_q^{M}(x) &:=  
\begin{dcases}
\dfrac{C}{\eta}\left(\dfrac{x-(\eta+(q-1)bh)}{\eta}\right)^2,\quad &x\in[(q-1)bh,\eta+(q-1)bh),\\
0,\quad &x\in[\eta+(q-1)bh,\eta+(qb-1)h],\\
\dfrac{C}{\eta}\left(\dfrac{x-(\eta+(qb-1)h)}{\eta}\right)^2,\quad &x\in(\eta+(qb-1)h,2\eta+(qb-1)h],
\end{dcases}\\
s_q^{M}(x) &:=  \left(1+\ii\dfrac{\sigma_q^{M}(x)}{\omega}\right)^{-1}, 
\end{align*}
and define $H_q^{M} \pmb v=\pmb g$ as the discrete problem of the continuous problem
\begin{align*}
\begin{dcases}
\left((s_q^{M}(x)\dfrac{\dd}{\dd x})^2+\dfrac{\omega^2}{c^2(x)}\right)v(x)=g(x),\quad & \forall x\in D_q^M,\\
v(x)=0, \quad & \forall x \in \partial D_q^M.
\end{dcases}
\end{align*}
We introduce the operator $\td{G}_q^{M}: \pmb y\to \pmb z$ as:
\begin{enumerate}
\item
  Introduce a vector $\pmb g$ defined on $X_q^{M}$ by setting $\pmb y$ to $X_q$
  and zero everywhere else.
\item
  Solve $H_q^{M} \pmb v=\pmb g$ on $X_q^{M}$. 
\item
  Set $\pmb z$ as the restriction of $\pmb v$ on $X_q$. 
\end{enumerate}
Then
\begin{align*}
  \td{\pmb u}_{q,q}:= \td{G}_q^{M} \pmb f_q. 
\end{align*}
Following the above discussion, the remaining components $\td{\pmb
  u}_{p,q}$ are defined recursively as
\begin{align*}
\td{\pmb u}_{p,q} &:=  \td{G}_{p}^{R} \td{\pmb u}_{p-1,q}^{R}, \quad \text{for }p=q+1,\dots,m,\\
\td{\pmb u}_{p,q} &:=  \td{G}_{p}^{L}\td{\pmb u}_{p+1,q}^{L},\quad \text{for } p=q-1,\dots,1. 
\end{align*}

\subsection{Accumulating the boundary values}
\begin{figure}[h!]
  \centering
  \subfigure[$\td {\pmb u}_p$ is a superposition of $\td {\pmb u}_{p,q},q=1,\dots,m$.]{            
    \label{fig:approxUseparate}
    \begin{tikzpicture}
      [x=2cm,y=2cm,>=latex]
      \draw[step=1](0,-0.1)grid(7,0.1);
      \draw(0.5,0)node[below=3]{$X_1$}(1.5,0)node[below=3]{$X_2$}(2.5,0)node[below=3]{$X_3$}(4.5,0)node[below=3]{$X_{m-2}$}(5.5,0)node[below=3]{$X_{m-1}$}(6.5,0)node[below=3]{$X_m$};
      \foreach \x in{0,...,5}
               {
                 \draw[->,gray,very thick](\x+1,2.5)--(\x+2,2.5);
                 \draw[<-,gray,very thick](\x,0.5)--(\x+1,0.5);
               }
               \draw[<->,gray,very thick](6,0.5)--(7,0.5);
               \draw[<->,gray,very thick](0,2.5)--(1,2.5);
               
               \foreach \x in{4,5,6}
               \draw[->,gray,very thick](\x,1.5)--(\x+1,1.5);
               \foreach \x in{0,1,2}
               \draw[<-,gray,very thick](\x,1.5)--(\x+1,1.5);
               \draw[<->,gray,very thick](3,1.5)--(4,1.5);
               
               \draw(6.5,0.5)node[fill=white]{$\pmb f_m$};
               \draw(0.5,2.5)node[fill=white]{$\pmb f_1$};
               
               \foreach \x in{0.5,1.5,2.5}
               \draw(3.5,\x)node[fill=white]{$\ldots$};
               
               \foreach \x in{0.5,1.5,2.5,4.5,5.5,6.5}
               \foreach \y in{1.0,2.0}
               \draw(\x,\y)node{$\vdots$};
               
               \draw(0.5,0.5)node[below=3]{$\td {\pmb u}_{1,m}$}(1.5,0.5)node[below=3]{$\td {\pmb u}_{2,m}$}(2.5,0.5)node[below=3]{$\td {\pmb u}_{3,m}$}(4.5,0.5)node[below=3]{$\td {\pmb u}_{m-2,m}$}(5.5,0.5)node[below=3]{$\td {\pmb u}_{m-1,m}$}(6.5,0.5)node[below=3]{$\td {\pmb u}_{m,m}$};
               
               \draw(3.5,0)node[fill=white]{$\ldots$};
               \draw(0.5,1.5)node[below=3]{$\td {\pmb u}_{1,q}$}(1.5,1.5)node[below=3]{$\td {\pmb u}_{2,q}$}(2.5,1.5)node[below=3]{$\td {\pmb u}_{3,q}$}(4.5,1.5)node[below=3]{$\td {\pmb u}_{m-2,q}$}(5.5,1.5)node[below=3]{$\td {\pmb u}_{m-1,q}$}(6.5,1.5)node[below=3]{$\td {\pmb u}_{m,q}$};
               
               \draw(3.5,0)node[fill=white]{$\ldots$};
               \draw(0.5,2.5)node[below=3]{$\td {\pmb u}_{1,1}$}(1.5,2.5)node[below=3]{$\td {\pmb u}_{2,1}$}(2.5,2.5)node[below=3]{$\td {\pmb u}_{3,1}$}(4.5,2.5)node[below=3]{$\td {\pmb u}_{m-2,1}$}(5.5,2.5)node[below=3]{$\td {\pmb u}_{m-1,1}$}(6.5,2.5)node[below=3]{$\td {\pmb u}_{m,1}$};
               
               \foreach \x in{0.5,1.5,2.5,4.5,5.5,6.5}
               \draw[->,>=latex,line width=3 pt](\x,2.7)--(\x,3.0);
               
               \draw(0.5,3.0)node[above]{$\td {\pmb u}_1$}(1.5,3.0)node[above]{$\td {\pmb u}_2$}(2.5,3.0)node[above]{$\td {\pmb u}_3$}(4.5,3.0)node[above]{$\td {\pmb u}_{m-2}$}(5.5,3.0)node[above]{$\td {\pmb u}_{m-1}$}(6.5,3.0)node[above]{$\td {\pmb u}_m$};
    \end{tikzpicture}
  }
  
  \subfigure[$\td {\pmb u}_p$ is a superposition of $\td {\pmb u}_{p,1:p-1}$, $\td {\pmb u}_{p,p}$ and $\td {\pmb u}_{p,p+1:m}$.]{
    \label{fig:approxUaccumulate}
    \begin{tikzpicture}
      [x=2cm,y=2cm,>=latex]
      \draw[step=1](0,-0.1)grid(7,0.1);
      \draw(0.5,0)node[below=3]{$X_1$}(1.5,0)node[below=3]{$X_2$}(2.5,0)node[below=3]{$X_3$}(4.5,0)node[below=3]{$X_{m-2}$}(5.5,0)node[below=3]{$X_{m-1}$}(6.5,0)node[below=3]{$X_m$};
      
      \foreach \x in{0,...,5}
               {
                 \draw[->,gray,very thick](\x+1,1.5)--(\x+2,1.5);
                 \draw[<-,gray,very thick](\x,0.5)--(\x+1,0.5);
               }
               \draw[<->,gray,very thick](6,0.5)--(7,0.5);
               \draw[<->,gray,very thick](0,1.5)--(1,1.5);
               
               \draw(6.5,0.5)node[fill=white]{$\pmb f_m$};
               \draw(0.5,1.5)node[fill=white]{$\pmb f_1$};
               
               \foreach \x in{0.5,1.5}
               \draw(3.5,\x)node[fill=white]{$\ldots$};
               
               \foreach \x in{1,...,5}
                        {
                          \draw[<->,gray,very thick](\x,0.5)--(\x+1,1.5);
                          \fill[white](3.35,0.7)rectangle(3.65,1.3);
                        }
                        
                        \draw(1.5,1.0)node[fill=white]{$\pmb f_2$}(2.5,1.0)node[fill=white]{$\pmb f_3$}(3.5,1.0)node[fill=white]{$\ldots$}(4.5,1.0)node[fill=white]{$\pmb f_{m-2}$}(5.5,1.0)node[fill=white]{$\pmb f_{m-1}$};
                        
                        \draw(1.5,1.0)node[below=3]{$\td {\pmb u}_{2,2}$}(2.5,1.0)node[below=3]{$\td {\pmb u}_{3,3}$}(4.5,1.0)node[below=3]{$\td {\pmb u}_{m-2,m-2}$}(5.5,1.0)node[below=3]{$\td {\pmb u}_{m-1,m-1}$};
                        
                        \draw(0.5,0.5)node[below=3]{$\td {\pmb u}_{1,2:m}$}(1.5,0.5)node[below=3]{$\td {\pmb u}_{2,3:m}$}(2.5,0.5)node[below=3]{$\td {\pmb u}_{3,4:m}$}(4.5,0.5)node[below=3]{$\td {\pmb u}_{m-2,m-1:m}$}(5.5,0.5)node[below=3]{$\td {\pmb u}_{m-1,m}$}(6.5,0.5)node[below=3]{$\td {\pmb u}_{m,m}$};
                        
                        \draw(3.5,0)node[fill=white]{$\ldots$};
                        \draw(0.5,1.5)node[below=3]{$\td {\pmb u}_{1,1}$}(1.5,1.5)node[below=3]{$\td {\pmb u}_{2,1}$}(2.5,1.5)node[below=3]{$\td {\pmb u}_{3,1:2}$}(4.5,1.5)node[below=3]{$\td {\pmb u}_{m-2,1:m-3}$}(5.5,1.5)node[below=3]{$\td {\pmb u}_{m-1,1:m-2}$}(6.5,1.5)node[below=3]{$\td {\pmb u}_{m,1:m-1}$};
                        
                        \foreach \x in{0.5,1.5,2.5,4.5,5.5,6.5}
                        \draw[->,>=latex,line width=3 pt](\x,1.7)--(\x,2.0);
                        
                        \draw(0.5,2.0)node[above]{$\td {\pmb u}_1$}(1.5,2.0)node[above]{$\td {\pmb u}_2$}(2.5,2.0)node[above]{$\td {\pmb u}_3$}(4.5,2.0)node[above]{$\td {\pmb u}_{m-2}$}(5.5,2.0)node[above]{$\td {\pmb u}_{m-1}$}(6.5,2.0)node[above]{$\td {\pmb u}_m$};
                        
    \end{tikzpicture}
  }
  \caption{This figure shows how the boundary values are accumulated after each step. The thin arrows indicate the transmission directions of the waves. The bold, up-pointing arrows symbolizes that summing up the corresponding waves on $X_p$ gives the superposition wave $\td {\pmb u}_p$.}
\end{figure}
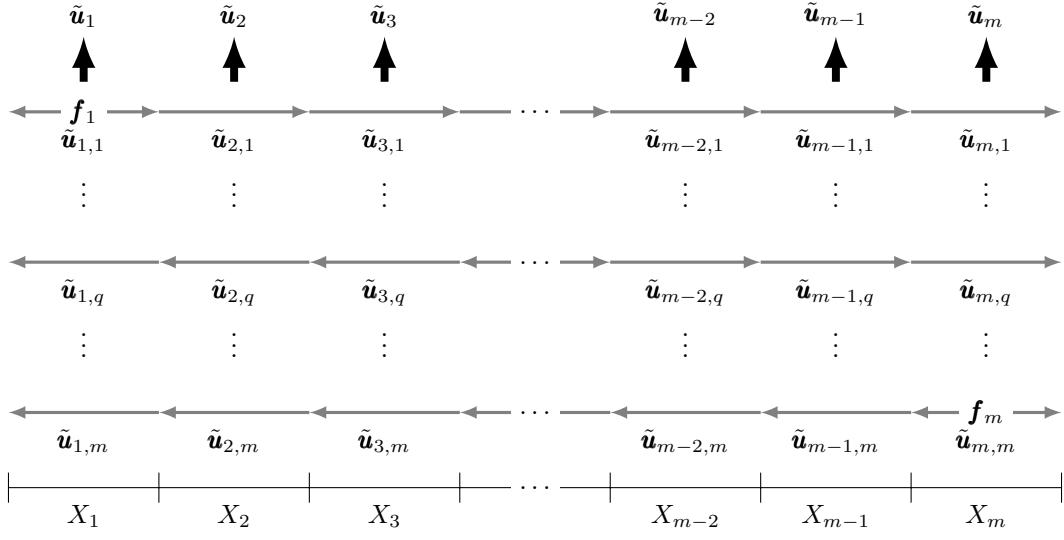
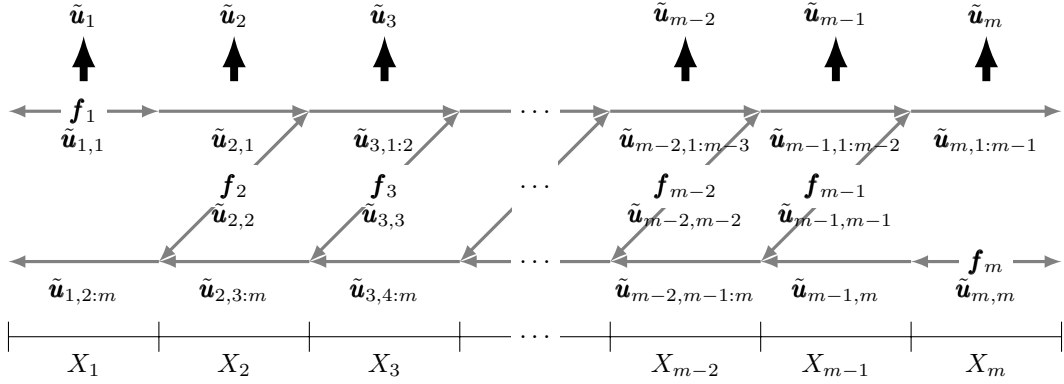

After all the above are done, an approximation of $\pmb u_p$ is
given by (see Figure \ref{fig:approxUseparate})
\begin{align*}
\td{\pmb u}_p:=\sum_{q=1}^m\td{\pmb u}_{p,q},\quad p=1,\dots,m.
\end{align*}

In the algorithm described above, the computation of each component
$\td{\pmb u}_{p,q}$ requires a separate solution of a problem of form
$H_p^{R} \pmb v =\pmb g$ or $H_p^{L} \pmb v =\pmb g$. Since there are $O(m^2)$ such
components, the algorithm is computationally expensive. A key
observation is that the computation associated with each $p$ can be
combined in one single shot by accumulating the boundary values of the waves.  More precisely, we define
\begin{align*}
\td{\pmb u}_{p,q_1:q_2}:=\sum_{t=q_1}^{q_2} \td{\pmb u}_{p,t}, 
\end{align*}
which is the total contribution of the waves generated by $\pmb
f_{q_1},\dots,\pmb f_{q_2}$ restricted to the grid $X_{p}$. The
quantity $\td{\pmb u}_{p,1:p-1}$, which is the total right-going wave
generated by $\pmb f_1,\dots,\pmb f_{p-1}$ upon $X_p$, can be computed
sequentially for $p=2,\dots,m$ without computing each component and
then adding them together as we described above, as long as we
accumulate the boundary values after each intermediate
step. Specifically, we first compute $\td{\pmb
  u}_{q,q}=\td{G}_q^{M}\pmb f_{q}$ for $q=1,\dots,m$.  This step is similar to what we did above. Then, to compute $\td{\pmb
  u}_{p,1:p-1}$ we carry out the following steps
\begin{align*}
  \td{\pmb u}_{p,1:p-1}=\td{G}_{p}^{R}\td{\pmb u}_{p-1,1:p-1}^{R},\quad \td{\pmb u}_{p,1:p}^{R}=\td{\pmb u}_{p,1:p-1}^{R}+\td{\pmb u}_{p,p}^{R}, \quad \text{for } p=2,\dots,m.  
\end{align*}
This means, before computing the total right-going wave $\td{\pmb u}_{p+1,1:p}$ on subdomain $X_{p+1}$, the boundary values of the previous right-going
waves, $\td{\pmb u}_{p,1:p-1}^{R}$ and $\td{\pmb u}_{p,p}^{R}$, are added together, so
that the the current right-going wave $\td{\pmb u}_{p+1,1:p}$ can be computed in one shot, eliminating the
trouble of solving the subproblems for many times and adding the
results together (see Figure \ref{fig:approxUaccumulate}).

For the left going waves $\td{\pmb u}_{p,p+1:m}$, a
similar process gives rise to the recursive formula
\begin{align*}
\td{\pmb u}_{p,p+1:m}=\td{G}_p^{L}\td{\pmb u}_{p+1,p+1:m}^{L},\quad \td{\pmb u}_{p,p:m}^{L}=\td{\pmb u}_{p,p}^{L}+\td{\pmb u}_{p,p+1:m}^{L}, \quad \text{for } p=m-1,\dots,1.
\end{align*}

Finally, each $\td{\pmb u}_p$ can be computed by summing $\td{\pmb
  u}_{p,1:p-1}$, $\td{\pmb u}_{p,p}$ and $\td{\pmb u}_{p,p+1:m}$
together (for the leftmost and the rightmost one, $\td{\pmb u}_1$ and
$\td{\pmb u}_m$, only two terms need to be summed), i.e.,
\begin{align*}
\td{\pmb u}_1 &= \td{\pmb u}_{1,1}+\td{\pmb u}_{1,2:m},\\
\td{\pmb u}_p &= \td{\pmb u}_{p,1:p-1}+\td{\pmb u}_{p,p}+\td{\pmb u}_{p,p+1:m}, \quad p=2,\dots,m-1,\\
\td{\pmb u}_m &= \td{\pmb u}_{m,1:m-1}+\td{\pmb u}_{m,m}. 
\end{align*}

We see that, by accumulating the boundary values after each intermediate step, we only need to solve $O(m)$ subproblems instead of $O(m^2)$.

In this algorithm, the approximation $\td{\pmb u}_{p}$ on each small
subdomain is divided into three parts. From a matrix point of view,
this is analogous to splitting the block matrix $G$ into its lower
triangular part, diagonal part and upper triangular part, and then
approximating each part as an operator to get the intermediate waves and then summing the intermediate results
together. This is why we call it the additive sweeping method. 

Equation \eqref{eqn:LDUsplit} shows an analogy of this procedure, where the matrix $G$ is split into $3m - 2$ blocks, each of which corresponds to a subproblem solving process:
\begin{align*}
\td{\pmb u}_{q,q} &\approx  {\pmb u}_{q,q} = G_{q,q} {\pmb f}_q,\quad q = 1,\dots, m,\\
\td{\pmb u}_{p,1:p-1} &\approx  {\pmb u}_{p,1:p-1} = \sum_{q=1}^{p-1} G_{p,q} {\pmb f}_q,\quad p = 2,\dots, m,\\
\td{\pmb u}_{p,p+1:m} &\approx  {\pmb u}_{p,p+1:m} = \sum_{q=p+1}^{m} G_{p,q} {\pmb f}_q,\quad p = 1,\dots, m-1.\\
\end{align*}

\begin{equation}
  \label{eqn:LDUsplit}
  \begin{bmatrix}
    \pmb u_1\\\pmb u_2\\ \ldots\\\pmb u_m
  \end{bmatrix}
  =
  \begin{bmatrix}
    \pmb u_{1,1}+\pmb u_{1,2:m}\\\pmb u_{2,1}+\pmb u_{2,2}+\pmb u_{2,3:m}\\ \ldots\\\pmb u_{m,1:m-1}+\pmb u_{m,m}
  \end{bmatrix}
  =
  \left[
    \begin{array}{ccccc}
      \multicolumn{1}{c|}{G_{1,1}} & G_{1,2} & \dots & & G_{1,m} \\ 
      \hline
      G_{2,1} & \multicolumn{1}{|c|}{G_{2,2}} & G_{2,3} & \ldots & G_{2,m} \\ 
      \hline
      & & \ldots & & \\ 
      \hline
      G_{m,1} & \ldots & • & G_{m,m-1} & \multicolumn{1}{|c}{G_{m,m}}
    \end{array} 
    \right]
  \begin{bmatrix}
    \pmb f_1\\\pmb f_2\\ \ldots\\\pmb f_m
  \end{bmatrix}
\end{equation}

When combined with standard iterative solvers, the approximation
algorithm serves as a preconditioner for Equation \eqref{eqn:1D} and it can be easily generalized to higher dimensions. In the following sections, we will discuss the details of the algorithm in 2D and 3D. To be structurally consistent,  we will keep the notations for 2D and 3D the same with the 1D case without causing ambiguity. Some of the key notations and concepts are listed below as a reminder to the reader:
\begin{itemize}
\item
$\{X_p\}_{p=1}^m$\quad The sliced partition of the discrete grid.
\item
$\{D_q^M\}_{q=1}^m$\quad The auxiliary domains with two-sided PML padding.
\item
$\{D_p^R\}_{p=2}^m$\quad The auxiliary domains with right-side PML padding.
\item
$\{D_p^L\}_{p=1}^{m-1}$\quad The auxiliary domains with left-side PML padding.
\item
$\{X_q^M\}_{q=1}^m$\quad $X_q$ with two-sided PML padding, the discretization of $D_q^M$.
\item
$\{X_p^R\}_{p=2}^m$\quad $X_p$ with right-side PML padding, the discretization of $D_p^R$.
\item
$\{X_p^L\}_{p=1}^{m-1}$\quad $X_p$ with left-side PML padding, the discretization of $D_p^L$.
\item
$\{\tilde{G}_q^M\}_{q=1}^m$ \quad The auxiliary Green's operators each of which maps the force on $X_q$ to the approximation of the wave field restricted to $X_q$.
\item
$\{\tilde{G}_p^R\}_{p=2}^m$ \quad The auxiliary Green's operators each of which maps the left boundary value to the approximated wave field restricted to $X_p$, which simulates the right-transmission of the waves.
\item
$\{\tilde{G}_p^L\}_{p=1}^{m-1}$ \quad The auxiliary Green's operators each of which maps the right boundary value to the approximated wave field restricted to $X_p$, which simulates the left-transmission of the waves.
\end{itemize}

\section{Preconditioner in 2D}
\label{sec:2D}

\subsection{Algorithm}

The domain of interest is $D=(0,1)^2$. We put PML on the two
opposite sides of the boundary, $x_2=0$ and $x_2=1$, to illustrate the
idea. The resulting equation is
\begin{align*}
  \begin{dcases}
    \left(\partial_1^2+(s(x_2)\partial_2)^2+\dfrac{\omega^2}{c^2(x)}\right)u(x)=f(x),&\quad \forall x=(x_1,x_2)\in D,\\
    u(x)=0,&\quad \forall x\in \partial D, 
  \end{dcases}
\end{align*}
We discretize $D$ with step size $h=1/(n+1)$ in each direction, which
results the Cartesian grid
\begin{align*}
  X:=\{(i_1h,i_2h):1\le i_1,i_2\le n\},  
\end{align*} 
and the discrete equation
\begin{equation}
  \label{eqn:2D}
  \begin{gathered}
    \dfrac{s_{i_2}}{h}\left(\dfrac{s_{i_2+1/2}}{h}(u_{i_1,i_2+1}-u_{i_1,i_2})-\dfrac{s_{i_2-1/2}}{h}(u_{i_1,i_2}-u_{i_1,i_2-1})\right)\\
    +\dfrac{u_{i_1+1,i_2}-2u_{i_1,i_2}+u_{i_1-1,i_2}}{h^2}+\dfrac{\omega^2}{c_{i_1,i_2}^2}u_{i_1,i_2}=f_{i_1,i_2}, \quad \forall 1\le i_1,i_2\le n,
  \end{gathered}
\end{equation}
where the subscript $(i_1,i_2)$ means that the corresponding function
is evaluated at $(i_1h,i_2h)$, and since $s(x_2)$ is a function of
$x_2$ only, we omit the $i_1$ subscript. $\pmb u$ and $\pmb f$ are
defined to be the column-major ordering of the discrete array $u$ and
$f$ on the grid $X$
\begin{align*}
  \pmb u:=[u_{1,1},\dots,u_{n,1},\dots,u_{n,n}]^T,\quad \pmb f:=[f_{1,1},\dots,f_{n,1},\dots,f_{n,n}]^T.
\end{align*}
Now \eqref{eqn:2D} can be written as $A\pmb u=\pmb f$.

We divide the grid into $m$ parts along the $x_2$ direction
\begin{align*}
  X_1 &:= \{(i_1h,i_2h):1\le i_1\le n,1 \le i_2 \le \gamma + b-1 \}, \\
  X_p &:= \{(i_1h,i_2h):1\le i_1\le n,\gamma + (p-1) b \le i_2 \le \gamma + pb-1 \},\quad p=2,\dots, m-1, \\
  X_m &:= \{(i_1h,i_2h):1\le i_1\le n,\gamma + (m-1) b \le i_2 \le 2 \gamma + mb-2 \}, 
\end{align*}
and we define $\pmb u_p$ and $\pmb f_p$ as the column-major ordering
restriction of $u$ and $f$ on $X_p$
\begin{align*}
  \pmb u_1 &:= [u_{1,1},\dots,u_{n,1},\dots,u_{n,\gamma+b-1}]^T,\\
  \pmb u_p &:= [u_{1,\gamma+(p-1)b},\dots,u_{n,\gamma+(p-1)b},\dots,u_{n,\gamma + pb-1}]^T, \quad p=2,\dots,m-1,\\
  \pmb u_m &:= [u_{1,\gamma + (m-1) b},\dots,u_{n,\gamma + (m-1) b},\dots,u_{n,2 \gamma + mb-2}]^T,\\
  \pmb f_1 &:= [f_{1,1},\dots,f_{n,1},\dots,f_{n,\gamma+b-1}]^T,\\
  \pmb f_p &:= [f_{1,\gamma+(p-1)b},\dots,f_{n,\gamma+(p-1)b},\dots,f_{n,\gamma + pb-1}]^T, \quad p=2,\dots,m-1,\\
  \pmb f_m &:= [f_{1,\gamma + (m-1) b},\dots,f_{n,\gamma + (m-1) b},\dots,f_{n,2 \gamma + mb-2}]^T,
\end{align*}
then $\pmb u=G\pmb f$ for $G=A^{-1}$ can be written as
\begin{align*}
\begin{bmatrix}
\pmb u_1\\\pmb u_2\\ \vdots\\\pmb u_m
\end{bmatrix}
=
\begin{bmatrix}
G_{1,1}&G_{1,2}&\ldots&G_{1,m}\\
G_{2,1}&G_{2,2}&\ldots&G_{2,m}\\
\vdots&\vdots&&\vdots\\
G_{m,1}&G_{m,2}&\ldots&G_{m,m}
\end{bmatrix}
\begin{bmatrix}
\pmb f_1\\\pmb f_2\\ \vdots\\\pmb f_m
\end{bmatrix}.
\end{align*}

\paragraph{Auxiliary domains.}
Following to the 1D case, the extended subdomains and the
corresponding left and right boundaries are defined by
\begin{align*}
  D_q^{M} &= (0,1)\times((q-1)bh,2\eta+(qb-1)h),\quad q=1,\dots,m,\\
  D_p^{R} &= (0,1)\times(\eta+((p-1)b-1)h,2\eta+(pb-1)h),\quad p=2,\dots,m,\\
  D_p^{L} &= (0,1)\times((p-1)bh,\eta+pbh),\quad p=1,\dots,m-1,\\
  \partial^{L} D_p^{R} &= (0,1)\times\{\eta+((p-1)b-1)h\},\quad p=2,\dots,m,\\
  \partial^{R} D_p^{L} &= (0,1)\times\{\eta+pbh\},\quad p=1,\dots,m-1.
\end{align*}
The extended grid for these domains are
\begin{align*}
  X_q^{M} &:= \{(i_1h,i_2h):1\le i_1\le n,(q-1)b+1\le i_2 \le 2\gamma+qb-1\}, \quad q=1,\dots,m, \\
  X_p^{R} &:= \{(i_1h,i_2h):1\le i_1\le n,\gamma+(p-1)b\le i_2\le 2\gamma+pb-2\}, \quad p=2,\dots,m, \\
  X_p^{L} &:= \{(i_1h,i_2h):1\le i_1\le n,(p-1)b+1\le i_2\le \gamma+pb-1\}, \quad p=1,\dots,m-1. 
\end{align*}

\paragraph{Auxiliary problems.}
For $q=1,\dots,m$, we define $H_q^{M} \pmb v=\pmb g$ to be the
discretization on $X_q^{M}$ of the problem
\begin{align*}
  &\begin{dcases}
     \left(\partial_1^2+(s_q^{M}(x_2)\partial_2)^2+\dfrac{\omega^2}{c^2(x)}\right)v(x)=g(x),\quad &\forall x\in D_q^{M},\\
     v(x)=0,\quad &\forall x\in \partial D_q^{M}. 
   \end{dcases}
\end{align*}
For $p=2,\dots,m$, $H_p^{R} \pmb v=\pmb g$ is the discretization on
$X_p^{R}$ of the problem
\begin{align*}
  &\begin{dcases}
     \left(\partial_1^2+(s_p^{R}(x_2)\partial_2)^2+\dfrac{\omega^2}{c^2(x)}\right)v(x)=0,\quad &\forall x\in D_p^{R},\\
     v(x)=w(x_1),\quad &\forall x\in \partial^{L} D_p^{R}, \\
     v(x)=0,\quad &\forall x\in \partial D_p^{R} \setminus \partial^{L} D_p^{R},
   \end{dcases}
\end{align*}
where $\pmb g:=(-1/h^2)[\pmb w^T,0,\dots,0]^T$ and $\pmb w:= [w_1,\dots,w_n]^T$ is the discrete value of $w(x_1)$. Finally, for
$p=1,\dots,m-1$, $H_p^{L} \pmb v=\pmb g$ is the discretization on
$X_p^{L}$ of the problem
\begin{align*}
  &\begin{dcases}
     \left(\partial_1^2+(s_p^{L}(x_2)\partial_2)^2+\dfrac{\omega^2}{c^2(x)}\right)v(x)=0,\quad &\forall x\in D_p^{L},\\
     v(x)=w(x_1),\quad &\forall x\in \partial^{R} D_p^{L}, \\
     v(x)=0,\quad &\forall x\in \partial D_p^{L} \setminus \partial^{R} D_p^{L},
   \end{dcases}
\end{align*}
where $\pmb g:=(-1/h^2)[0,\dots,0,\pmb w^T]^T$ and $\pmb w:= [w_1,\dots,w_n]^T$.

\paragraph{Auxiliary Green's operators.}
For $q=1,\dots,m$, we define $\td{G}_q^{M}:\pmb y\mapsto \pmb z$ to be
the operator defined by the following operations:
\begin{enumerate}
\item
  Introduce a vector $\pmb g$ defined on $X_q^{M}$ by setting $\pmb y$ to $X_q$ and zero everywhere else. 
\item
  Solve $H_q^{M} \pmb v=\pmb g$ on $X_q^{M}$. 
\item
  Set $\pmb z$ as the restriction of $\pmb v$ on $X_q$. 
\end{enumerate}
For $p=2,\dots,m$, the operators $\td{G}_p^{R}:\pmb w\mapsto \pmb z$
is given by:
\begin{enumerate}
\item
  Set $\pmb g=(-1/h^2)[\pmb w^T,0,\dots,0]^T$. 
\item
  Solve $H_p^{R} \pmb v =\pmb g$ on $X_p^{R}$. 
\item
  Set $\pmb z$ as the restriction of $\pmb v$ on $X_p$. 
\end{enumerate}
Finally, for $p=1,\dots,m-1$, $\td{G}_p^{L}:\pmb w\mapsto \pmb z$ is
defined as:
\begin{enumerate}
\item
  Set $\pmb g=(-1/h^2)[0,\dots,0,\pmb w^T]^T$. 
\item
  Solve $H_p^{L} \pmb v =\pmb g$ on $X_p^{L}$. 
\item
  Set $\pmb z$ as the restriction of $\pmb v$ on $X_p$. 
\end{enumerate}

\paragraph{Putting together.}
Similar to the previous section, we introduce the left boundary value
$\pmb g^{L}$ and the right boundary value $\pmb g^{R}$ for a
column-major ordering array
$\pmb g=[g_{1,1},\dots,g_{s_1,1},\dots,g_{s_1,s_2}]^T$
induced from some grid with size $s_1\times s_2$ by
\begin{align*}
  \pmb g^{L}:=[g_{1,1},\dots,g_{s_1,1}]^T, \quad \pmb g^{R}:=[g_{1,s_2},\dots,g_{s_1,s_2}]^T. \\
\end{align*}
Then the approximations for $\pmb u_p,p=1,\dots,m$, can be defined step by step as
\begin{align*}
  \td{\pmb u}_{q,q} &:= \td{G}_q^{M} \pmb f_q,\quad q=1,\dots,m,\\
  \td{\pmb u}_{p,1:p-1} &:= \td{G}_{p}^{R}\td{\pmb u}_{p-1,1:p-1}^{R},\quad \td{\pmb u}_{p,1:p}^{R}:=\td{\pmb u}_{p,1:p-1}^{R}+\td{\pmb u}_{p,p}^{R},\quad \text{for } p=2,\dots,m,\\
  \td{\pmb u}_{p,p+1:m} &:= \td{G}_p^{L}\td{\pmb u}_{p+1,p+1:m}^{L},\quad \td{\pmb u}_{p,p:m}^{L}:=\td{\pmb u}_{p,p}^{L}+\td{\pmb u}_{p,p+1:m}^{L}, \quad \text{for } p=m-1,\dots,1,\\
  \td{\pmb u}_1 &:= \td{\pmb u}_{1,1}+\td{\pmb u}_{1,2:m},\\
  \td{\pmb u}_p &:= \td{\pmb u}_{p,1:p-1}+\td{\pmb u}_{p,p}+\td{\pmb u}_{p,p+1:m},\quad p=2,\dots,m-1,\\
  \td{\pmb u}_m &:= \td{\pmb u}_{m,1:m-1}+\td{\pmb u}_{m,m}. 
\end{align*}

To solve the subproblems $H_q^{M} \pmb v=\pmb g$, $H_p^{R}\pmb v=\pmb
g$ and $H_p^{L}\pmb v=\pmb g$, we notice that they are indeed quasi-1D
problems since $\gamma$ and $b$ are some small constants. Therefore,
for each one of them, we can reorder the system by grouping
the elements along dimension 2 first and then dimension 1, which
results a banded linear system that can be solved by the LU
factorization efficiently. These factorization processes induce the
factorizations for the operators $\td{G}_q^{M}$, $\td{G}_p^{R}$ and
$\td{G}_p^{L}$ symbolically, which leads to our setup algorithm of the
preconditioner in 2D as described in Algorithm \ref{alg:2dsetup} and
the application algorithm as described in Algorithm \ref{alg:2dapp}.

\begin{algorithm}[h!]
  \caption{Construction of the 2D additive sweeping preconditioner of
    the Equation \eqref{eqn:2D}. Complexity
    $=O(n^2(b+\gamma)^3/b)=O(N(b+\gamma)^3/b)$.}
  \label{alg:2dsetup}
  \begin{algorithmic}
    \FOR {$q=1,\dots,m$}
    \STATE Construct the LU factorization of $H_q^{M}$, which defines $\td{G}_q^{M}$.
    \ENDFOR
    \FOR {$p=2,\dots,m$}
    \STATE Construct the LU factorization of $H_p^{R}$, which defines $\td{G}_p^{R}$.
    \ENDFOR
    \FOR {$p=1,\dots,m-1$}
    \STATE Construct the LU factorization of $H_p^{L}$, which defines $\td{G}_p^{L}$.
    \ENDFOR
  \end{algorithmic}
\end{algorithm}

\begin{algorithm}[h!]
  \caption{Computation of $\td{\pmb u}\approx G \pmb f$ using the
    preconditioner from Algorithm \ref{alg:2dsetup}. Complexity
    $=O(n^2(b+\gamma)^2/b)=O(N(b+\gamma)^2/b)$.}
  \label{alg:2dapp}
  \begin{algorithmic}
    \FOR {$q=1,\dots,m$}
    \STATE
    $\td{\pmb u}_{q,q}=\td{G}_q^{M} \pmb f_q$
    \ENDFOR
    \FOR {$p=2,\dots,m$}
    \STATE
    $\td{\pmb u}_{p,1:p-1}=\td{G}_{p}^{R}\td{\pmb u}_{p-1,1:p-1}^{R}$\\
    $\td{\pmb u}_{p,1:p}^{R}=\td{\pmb u}_{p,1:p-1}^{R}+\td{\pmb u}_{p,p}^{R}$
    \ENDFOR
    \FOR {$p=m-1,\dots,1$}
    \STATE
    $\td{\pmb u}_{p,p+1:m}=\td{G}_p^{L}\td{\pmb u}_{p+1,p+1:m}^{L}$\\
    $\td{\pmb u}_{p,p:m}^{L}=\td{\pmb u}_{p,p}^{L}+\td{\pmb u}_{p,p+1:m}^{L}$
    \ENDFOR
    \STATE
    $\td{\pmb u}_1=\td{\pmb u}_{1,1}+\td{\pmb u}_{1,2:m}$\\
    \FOR {$p=2,\dots,m-1$}
    \STATE
    $\td{\pmb u}_p=\td{\pmb u}_{p,1:p-1}+\td{\pmb u}_{p,p}+\td{\pmb u}_{p,p+1:m}$
    \ENDFOR
    \STATE
    $\td{\pmb u}_m=\td{\pmb u}_{m,1:m-1}+\td{\pmb u}_{m,m}$ 
  \end{algorithmic}
\end{algorithm}

To analyze the complexity, we note that, in the setup process, there
are $O(n/b)$ subproblems, each of which is a quasi-1D problem with
$O(\gamma+b)$ layers along the second dimension. Therefore, the setup
cost of each subproblem by the LU factorization is $O(n(\gamma+b)^3)$
and the application cost is $O(n(\gamma+b)^2)$. So the total setup
cost is $O(n^2(\gamma+b)^3/b)$. Besides, one needs to solve each
subproblem once during the application process so the total
application cost is $O(n^2(\gamma+b)^2/b)$.

There are some differences when implementing the method practically:
\begin{enumerate}
\item
  In the above setting, PMLs are put only on two opposite sides of the
  unit square for illustration purpose. In reality, PMLs can be put on
  other sides of the domain if needed. As long as there are two
  opposite sides with PML boundary condition, the method can be
  implemented.
\item
  The thickness of the auxiliary PMLs introduced in the interior part
  of the domain needs not to be the same with the thickness of the PML
  at the boundary. In fact, the thickness of the auxiliary PML is
  typically thinner in order to improve efficiency.
\item
  The widths of the subdomains are completely arbitrary and they need
  not to be the same. Practically, the widths can be chosen to be
  larger for subdomains where the velocity field varies heavily.
\item
  The symmetric version of the equation can be adopted to save memory
  and computational cost.
\end{enumerate}

\subsection{Numerical results}
\label{sec:2Dnumerical}

\begin{figure}[h!]
  \centering
  \subfigure[]
  {\includegraphics[width=0.32\textwidth]{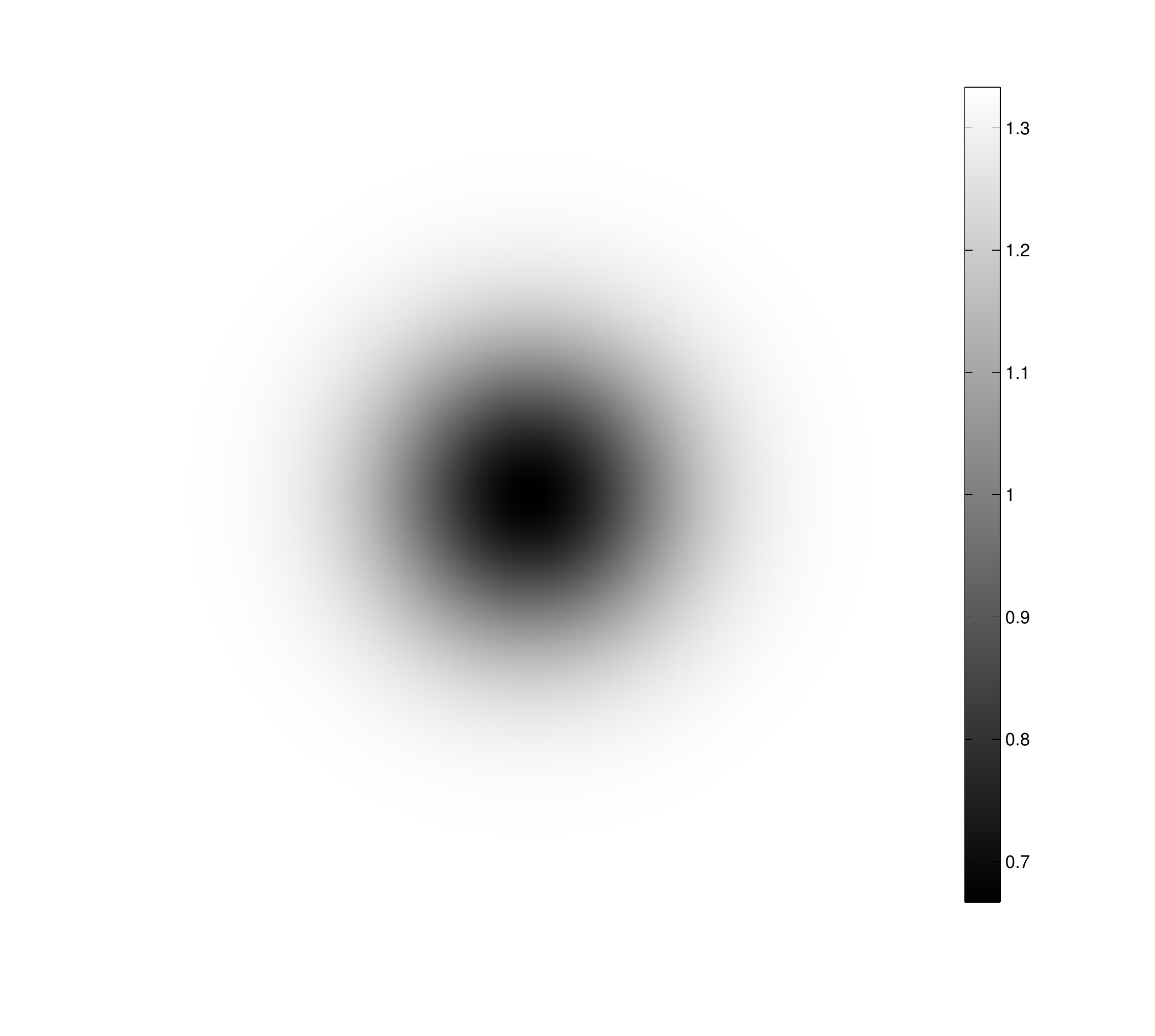}}
  \subfigure[]
  {\includegraphics[width=0.32\textwidth]{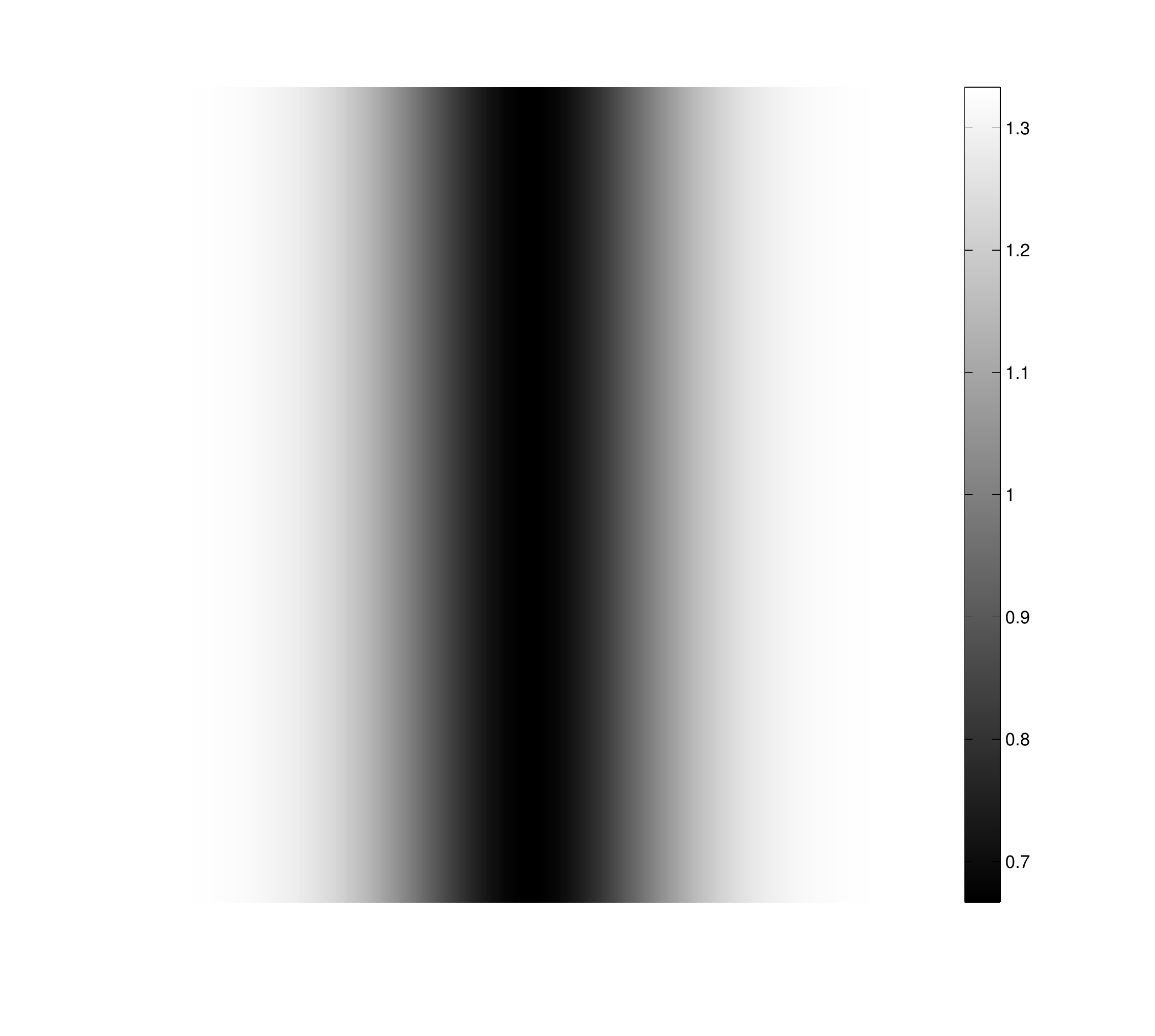}}
  \subfigure[]
  {\includegraphics[width=0.32\textwidth]{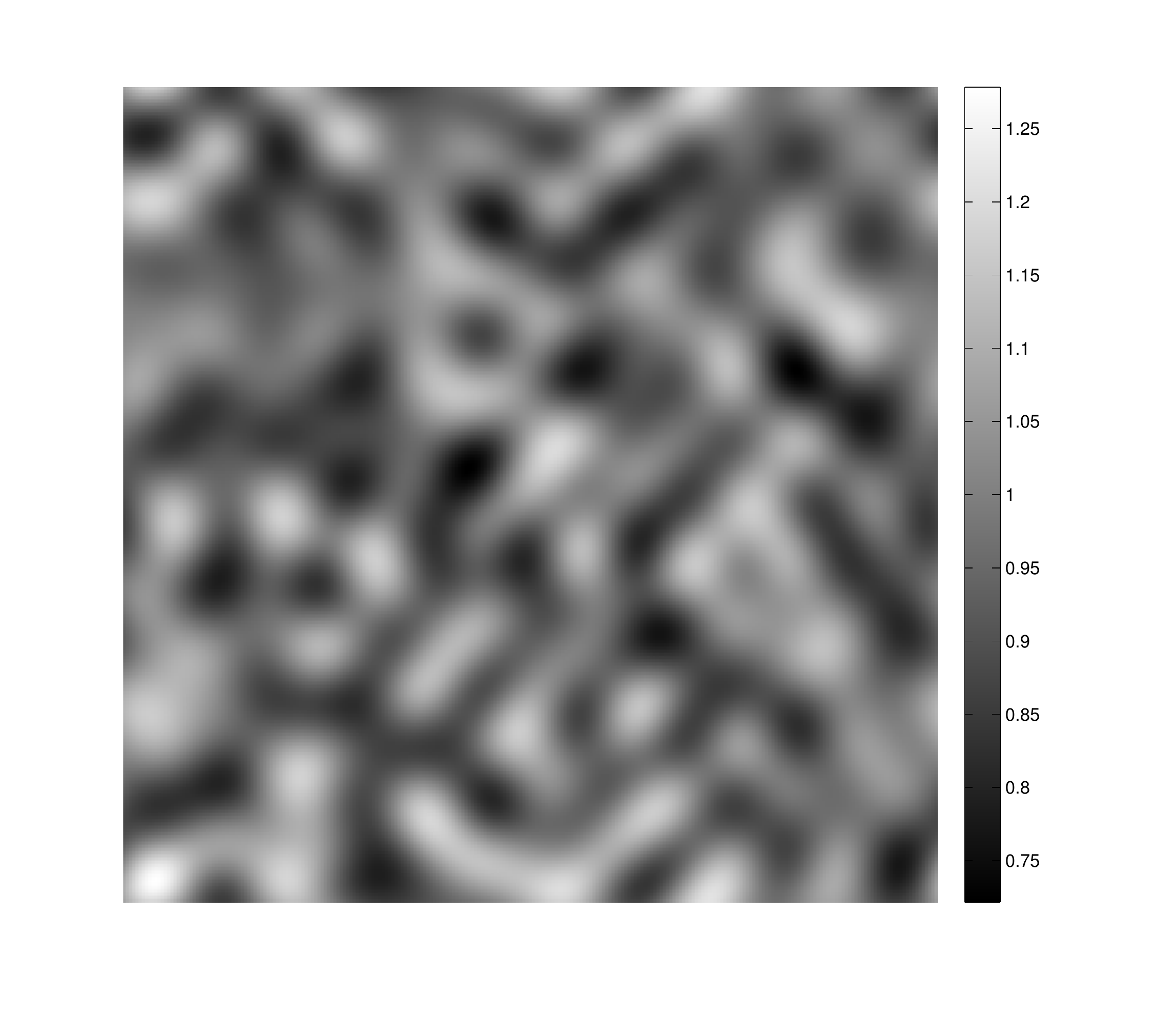}}
  \caption{The three velocity fields tested in 2D.}
  \label{fig:2D}
\end{figure}

Here, we present some numerical results in 2D to illustrate the
efficiency of the algorithm. The proposed method is implemented in
MATLAB and the tests are performed on a 2.0 GHz computer with 256 GB
memory. GMRES is used as the iterative solver with relative residual
equal to $10^{-3}$ and restart value equal to $40$. PMLs are put on
all sides of the unit square. The velocity fields tested are given in
Figure \ref{fig:2D}:
\begin{enumerate}[(a)]
\item
  A converging lens with a Gaussian profile at the center of the domain.
\item
  A vertical waveguide with a Gaussian cross-section.
\item
  A random velocity field.
\end{enumerate}

For each velocity field, two external forces are tested:
\begin{enumerate}[(a)]
\item
  A Gaussian point source centered at $(1/2,1/8)$.
\item
  A Gaussian wave packet with wavelength comparable to the typical
  wavelength of the domain. The packet centers at $(1/8,1/8)$ and
  points to the direction $(1/\sqrt{2},1/\sqrt{2})$.
\end{enumerate}

\begin{table}[h!]
  \centering
  \begin{overpic}
    [width=0.45\textwidth]{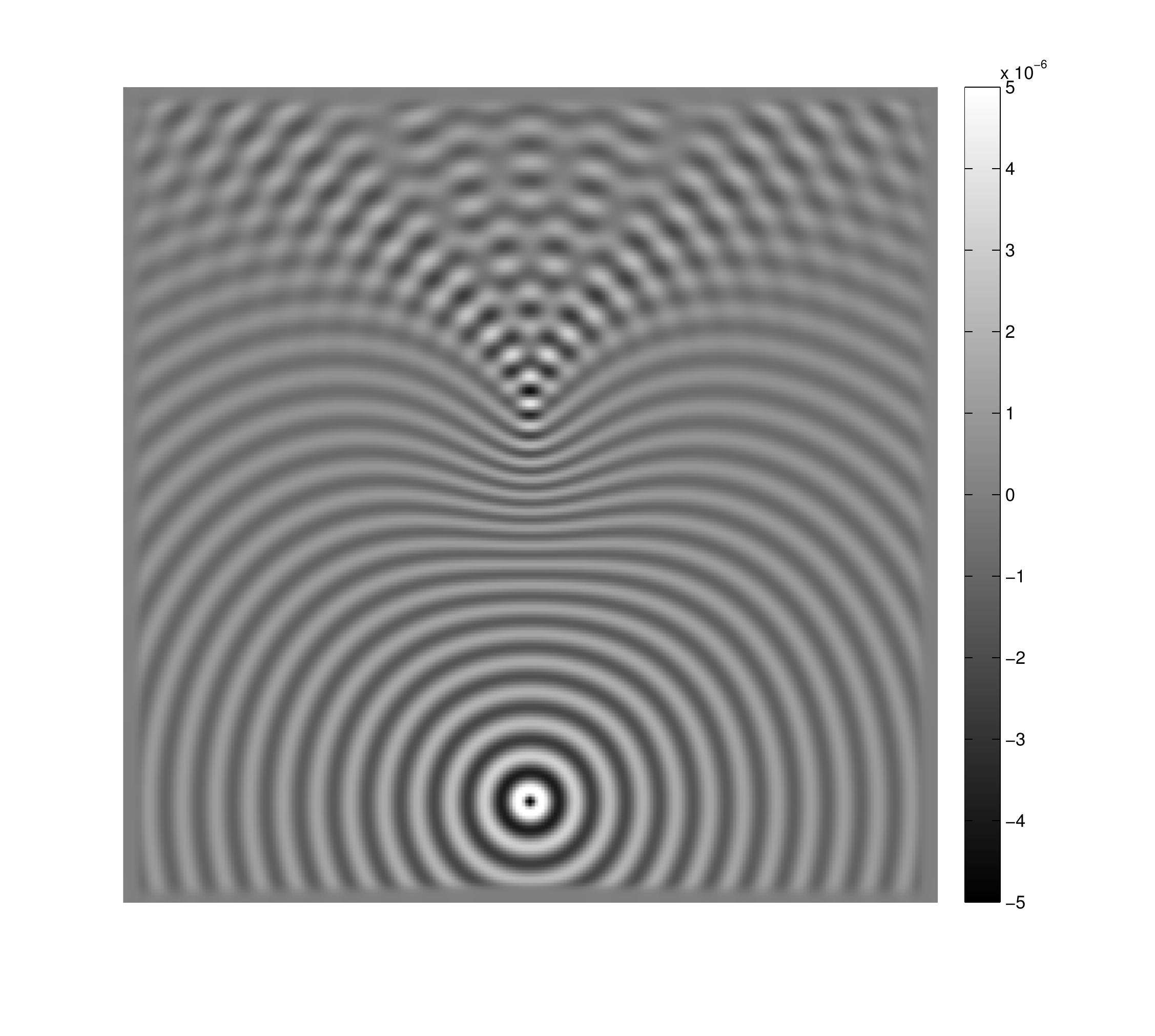}
    \put(40,6){force (a)}
  \end{overpic}
  \begin{overpic}
    [width=0.45\textwidth]{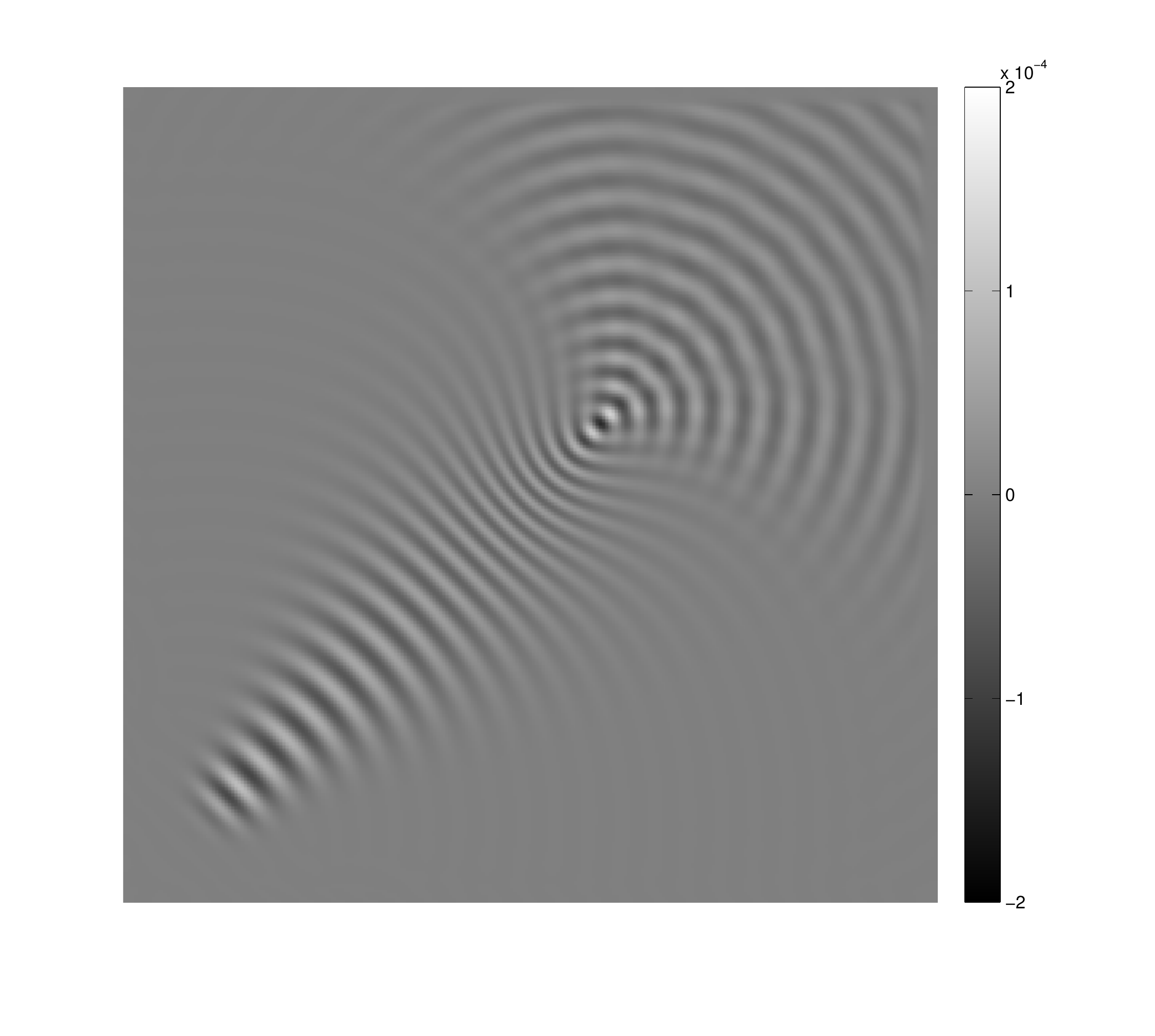}
    \put(40,6){force (b)}
  \end{overpic}
  \begin{tabular}{ccc|cc|cc}
    \hline 
    \multicolumn{3}{c|}{velocity field (a)} & \multicolumn{2}{c|}{force (a)} & \multicolumn{2}{c}{force (b)} \\ 
    \hline 
    $\omega/(2\pi)$ & $N$ & $T_{\text{setup}}$ & $N_\text{iter}$ & $T_\text{solve}$ & $N_\text{iter}$ & $T_\text{solve}$ \\ 
    \hline 
    16 & $127^2$ & 8.1669e$-$01 &  4 & 5.3199e$-$01 &  4 & 2.5647e$-$01 \\
    32 & $255^2$ & 3.4570e$+$00 &4& 7.3428e$-$01 &4&  7.2807e$-$01\\ 
    64 & $511^2$ & 1.5150e$+$01  &5&  3.6698e$+$00  &4& 3.7239e$+$00\\ 
    128 & $1023^2$ & 6.2713e$+$01 &5&  1.6812e$+$01 &4&  1.6430e$+$01\\
    256 & $2047^2$ & 2.6504e$+$02 &6&   7.8148e$+$01 &4& 5.6936e$+$01\\
    \hline 
  \end{tabular} 
  \caption{Results for velocity field (a) in 2D. Solutions with $\omega/(2\pi)=32$ are presented.}
  \label{tab:2domegaa}
\end{table}
\begin{table}[h!]
  \centering
  \begin{overpic}
    [width=0.45\textwidth]{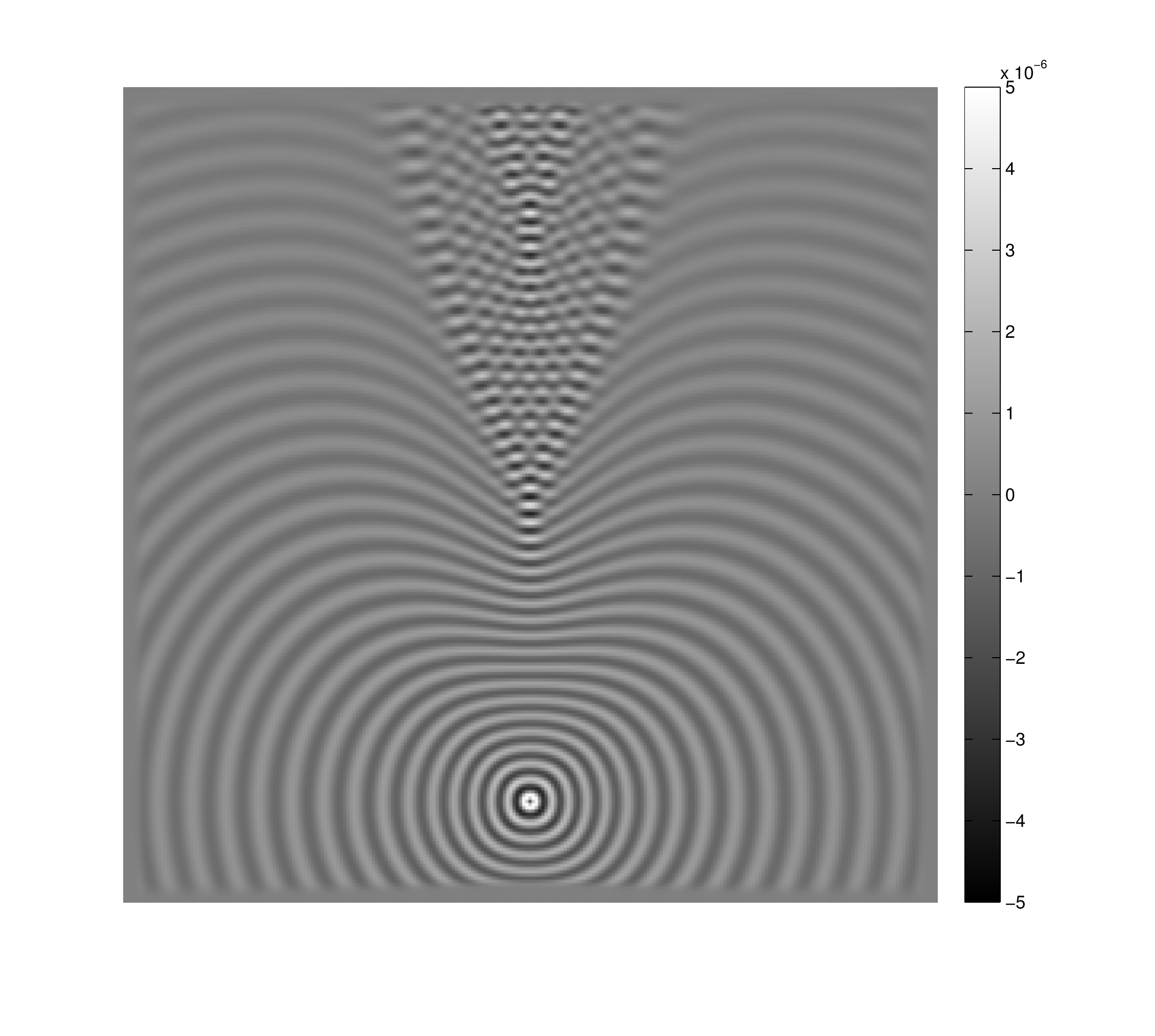}
    \put(40,6){force (a)}
  \end{overpic}
  \begin{overpic}
    [width=0.45\textwidth]{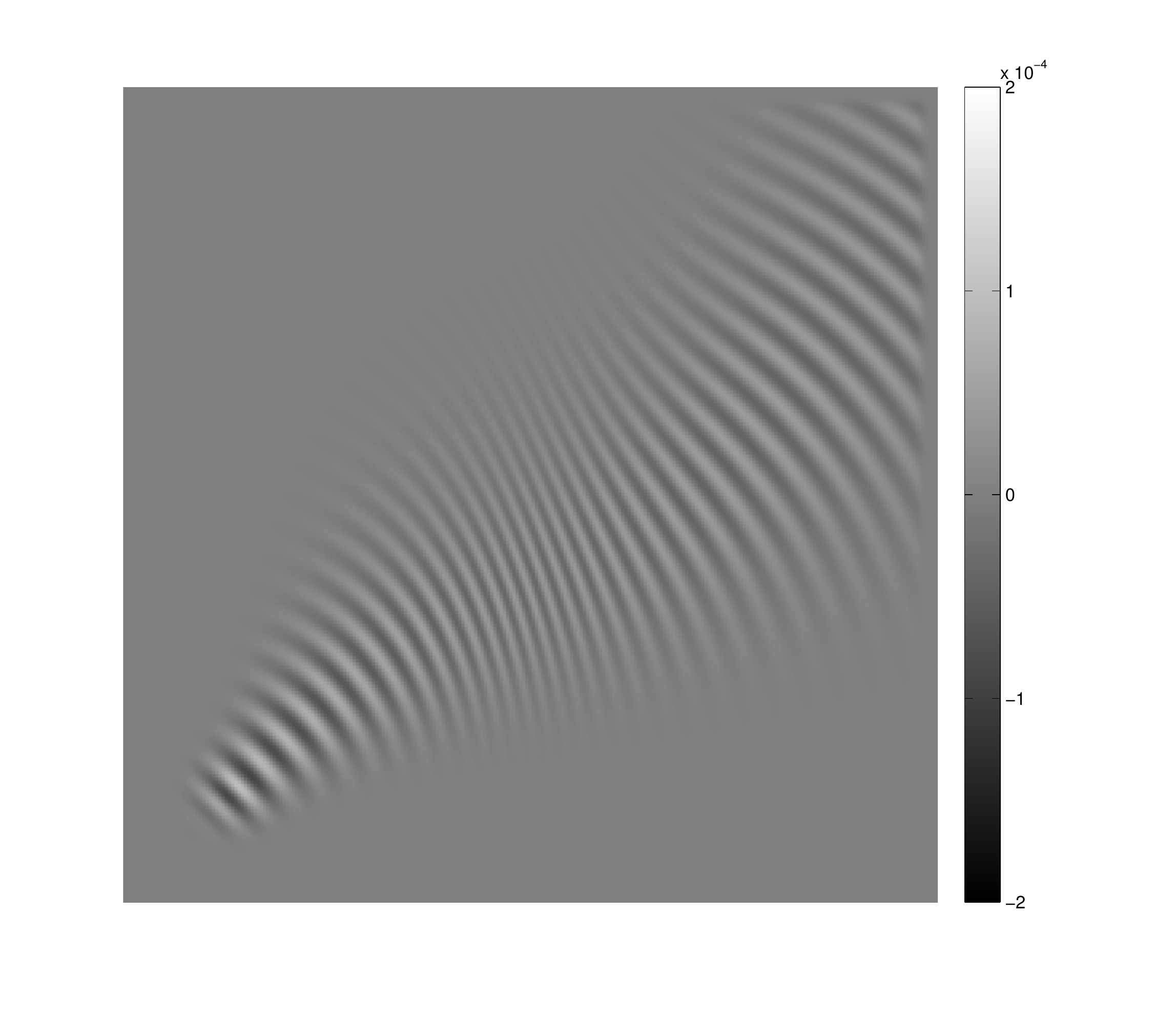}
    \put(40,6){force (b)}
  \end{overpic}
  \begin{tabular}{ccc|cc|cc}
    \hline 
    \multicolumn{3}{c|}{velocity field (b)} & \multicolumn{2}{c|}{force (a)} & \multicolumn{2}{c}{force (b)} \\ 
    \hline 
    $\omega/(2\pi)$ & $N$ & $T_{\text{setup}}$ & $N_\text{iter}$ & $T_\text{solve}$ & $N_\text{iter}$ & $T_\text{solve}$ \\ 
    \hline 
    16 & $127^2$ & 7.0834e$-$01 &6& 2.9189e$-$01 &4& 1.9408e$-$01\\
    32 & $255^2$ & 3.2047e$+$00 &8& 1.6147e$+$00  &4&  7.9303e$-$01\\ 
    64 & $511^2$ & 1.4079e$+$01 &8& 6.3057e$+$00   &4&  3.9008e$+$00\\ 
    128 & $1023^2$ & 6.0951e$+$01 &8&  2.9097e$+$01  &4&  1.5287e$+$01\\
    256 & $2047^2$ & 2.6025e$+$02  &8&  1.1105e$+$02  &5& 7.2544e$+$01\\
    \hline 
  \end{tabular} 
  \caption{Results for velocity field (b) in 2D. Solutions with $\omega/(2\pi)=32$ are presented.}
  \label{tab:2domegab}
\end{table}
\begin{table}[h!]
\centering
\begin{overpic}
  [width=0.45\textwidth]{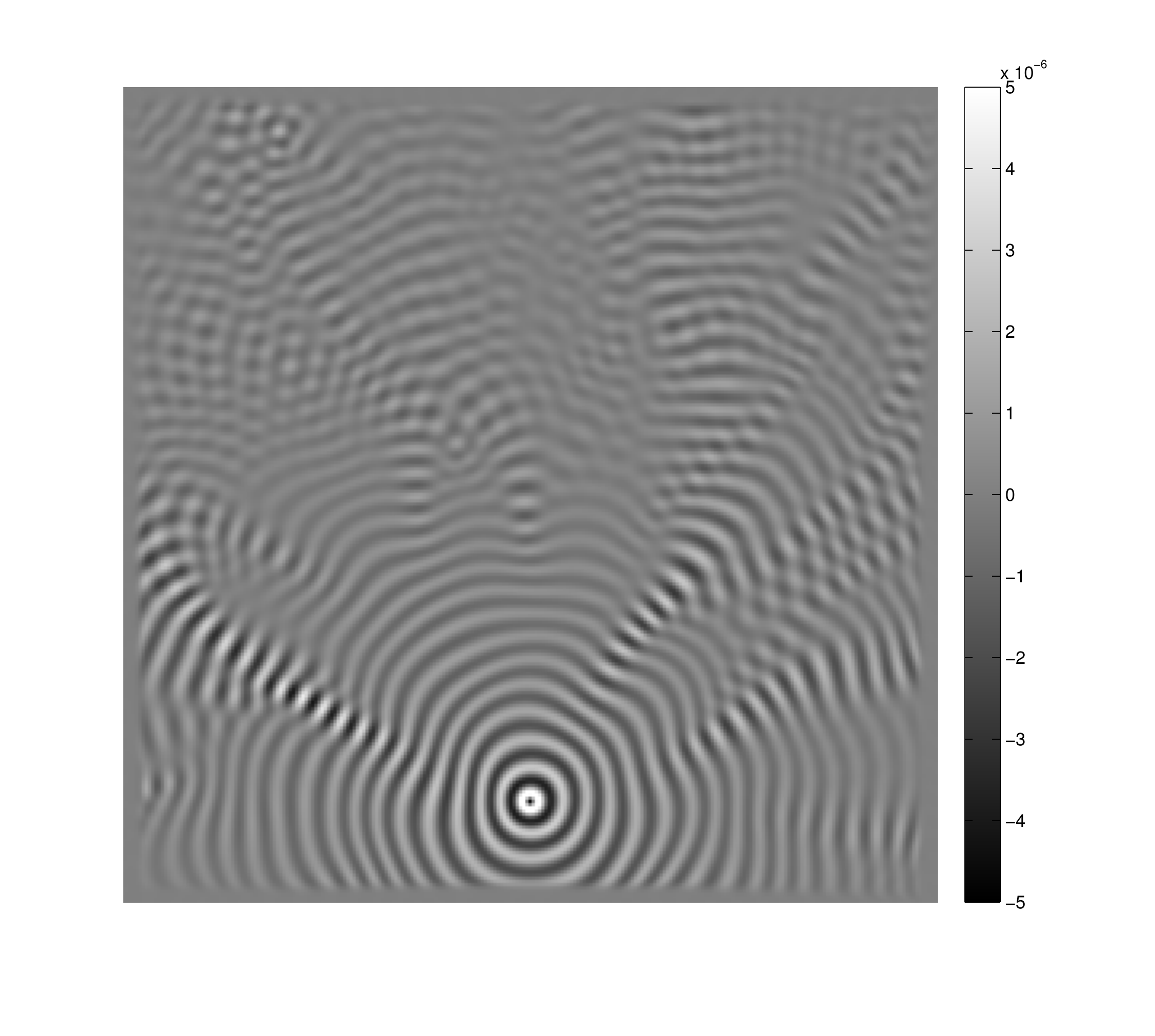}
  \put(40,6){force (a)}
\end{overpic}
\begin{overpic}
  [width=0.45\textwidth]{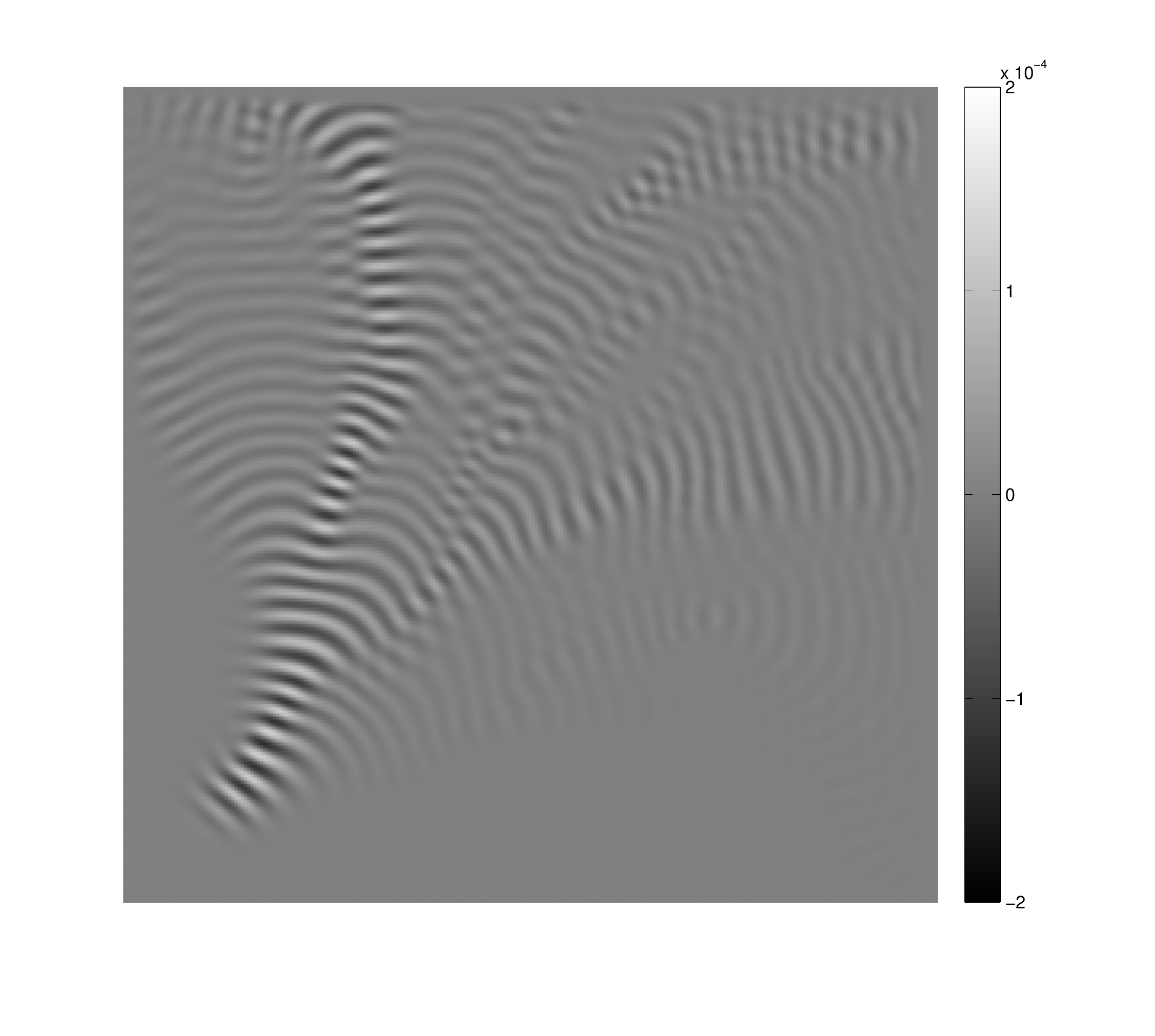}
  \put(40,6){force (b)}
\end{overpic}
\begin{tabular}{ccc|cc|cc}
  \hline 
  \multicolumn{3}{c|}{velocity field (c)} & \multicolumn{2}{c|}{force (a)} & \multicolumn{2}{c}{force (b)} \\ 
  \hline 
  $\omega/(2\pi)$ & $N$ & $T_{\text{setup}}$ & $N_\text{iter}$ & $T_\text{solve}$ & $N_\text{iter}$ & $T_\text{solve}$ \\ 
  \hline 
  16 & $127^2$ & 7.0495e$-$01  &5&  2.4058e$-$01  &6& 2.8347e$-$01\\
  32 & $255^2$ & 3.1760e$+$00 &5&  1.0506e$+$00  &5&  9.9551e$-$01\\ 
  64 & $511^2$ & 1.4041e$+$01   &6&  4.7083e$+$00  &7&  6.7852e$+$00\\
  128 & $1023^2$ & 6.1217e$+$01  &6&  1.8652e$+$01  &6&   1.9792e$+$01\\
  256 & $2047^2$ & 2.5762e$+$02  &8&  1.1214e$+$02  &6&  8.6936e$+$01\\
  \hline 
\end{tabular} 
\caption{Results for velocity field (c) in 2D. Solutions with $\omega/(2\pi)=32$ are presented.}
\label{tab:2domegac}
\end{table}

In these tests, each typical wavelength is discretized with 8
points. The width of the PML at the boundary and the one of the PMLs introduced
in the interior parts of the domain are both $9h$, i.e., $\gamma=9$. The number of layers in each interior subdomain is $b=8$,
the number of layers in the leftmost subdomain is $b+\gamma-1=16$ and
the one in the rightmost is $b+\gamma-2=15$. 

We vary the typical wave number $\omega/(2\pi)$ and test the behavior
of the algorithm. The test results are presented in Tables
\ref{tab:2domegaa}, \ref{tab:2domegab} and
\ref{tab:2domegac}. $T_\text{setup}$ is the setup time of the
algorithm in seconds. $T_{\text{solve}}$ is the total solve time in
seconds and $N_{\text{iter}}$ is the iteration number. From these
tests we see that the setup time scales like $O(N)$ as well as the
solve time per iteration, which is consistent with the algorithm
complexity analysis. The iteration number remains constant or grows at
most logarithmically, which shows the efficiency of the
preconditioner.

\section{Preconditioner in 3D}
\label{sec:3D}

\subsection{Algorithm}

In this section we briefly state the preconditioner in 3D case. The
domain of interest is $D=(0,1)^3$. PMLs are put on two opposite faces
of the unit cube, $x_3=0$ and $x_3=1$, which results the equation
\begin{align*}
  \begin{dcases}
    \left(\partial_1^2+\partial_2^2+(s(x_3)\partial_3)^2+\dfrac{\omega^2}{c^2(x)}\right)u(x)=f(x),&\quad \forall x=(x_1,x_2,x_3)\in D,\\
    u(x)=0,&\quad \forall x\in \partial D, 
  \end{dcases}
\end{align*}
Discretizing $D$ with step size $h=1/(n+1)$ gives the grid
\begin{align*}
  X:=\{(i_1h,i_2h,i_3h):1\le i_1,i_2,i_3\le n\},  
\end{align*} 
and the discrete equation
\begin{equation}
  \label{eqn:3D}
  \begin{gathered}
    \dfrac{s_{i_3}}{h}\left(\dfrac{s_{i_3+1/2}}{h}(u_{i_1,i_2,i_3+1}-u_{i_1,i_2,i_3})-\dfrac{s_{i_3-1/2}}{h}(u_{i_1,i_2,i_3}-u_{i_1,i_2,i_3-1})\right)\\
    +\dfrac{u_{i_1+1,i_2,i_3}-2u_{i_1,i_2,i_3}+u_{i_1-1,i_2,i_3}}{h^2}+\dfrac{u_{i_1,i_2+1,i_3}-2u_{i_1,i_2,i_3}+u_{i_1,i_2-1,i_3}}{h^2}\\
    +\dfrac{\omega^2}{c_{i_1,i_2,i_3}^2}u_{i_1,i_2,i_3}=f_{i_1,i_2,i_3}, \quad \forall 1\le i_1,i_2\le n. 
  \end{gathered}
\end{equation}
$\pmb u$ and $\pmb f$ are defined as the column-major ordering of $u$ and $f$ on the grid $X$
\begin{align*}
  \pmb u:=[u_{1,1,1},\dots,u_{n,1,1},\dots,u_{n,n,1},\dots,u_{n,n,n}]^T,\quad \pmb f:=[f_{1,1,1},\dots,f_{n,1,1},\dots,f_{n,n,1},\dots,f_{n,n,n}]. 
\end{align*}
$X$ is divided into $m$ parts along the $x_3$ direction
\begin{align*}
  X_1 &:= \{(i_1h,i_2h,i_3h):1\le i_1\le n,1\le i_2\le n,1 \le i_3 \le \gamma + b-1 \}, \\
  X_p &:= \{(i_1h,i_2h,i_3h):1\le i_1\le n,1\le i_2\le n,\gamma + (p-1) b \le i_3 \le \gamma + pb-1 \}, \quad p=2,\dots, m-1, \\
  X_m &:= \{(i_1h,i_2h,i_3h):1\le i_1\le n,1\le i_2\le n,\gamma + (m-1) b \le i_3 \le 2 \gamma + mb-2 \}. 
\end{align*}
$\pmb u_p$ and $\pmb f_p$ are the column-major ordering restrictions of  $u$ and $f$ on $X_p$ 
\begin{align*}
\pmb u_1 &:= [u_{1,1,1},\dots,u_{n,1,1},\dots,u_{n,n,1},\dots,u_{n,n,\gamma+b-1}]^T,\\
\pmb u_p &:= [u_{1,1,\gamma+(p-1)b},\dots,u_{n,1,\gamma+(p-1)b},\dots,u_{n,n,\gamma+(p-1)b},\dots,u_{n,n,\gamma + pb-1}]^T, \quad p=2,\dots,m-1,\\
\pmb u_m &:= [u_{1,1,\gamma + (m-1) b},\dots,u_{n,1,\gamma + (m-1) b},\dots,u_{n,n,\gamma + (m-1) b},\dots,u_{n,n,2 \gamma + mb-2}]^T,\\
\pmb f_1 &:= [f_{1,1,1},\dots,f_{n,1,1},\dots,f_{n,n,1},\dots,f_{n,n,\gamma+b-1}]^T,\\
\pmb f_p &:= [f_{1,1,\gamma+(p-1)b},\dots,f_{n,1,\gamma+(p-1)b},\dots,f_{n,n,\gamma+(p-1)b},\dots,f_{n,n,\gamma + pb-1}]^T, \quad p=2,\dots,m-1,\\
\pmb f_m &:= [f_{1,1,\gamma + (m-1) b},\dots,f_{n,1,\gamma + (m-1) b},\dots,f_{n,n,\gamma + (m-1) b},\dots,f_{n,n,2 \gamma + mb-2}]^T. 
\end{align*}

\paragraph{Auxiliary domains.}
The extended subdomains, the extended grids, and the corresponding
left and right boundaries are defined by
\begin{align*}
D_q^{M} &:= (0,1)\times(0,1)\times((q-1)bh,2\eta+(qb-1)h), \quad q=1,\dots,m,\\
D_p^{R} &:= (0,1)\times(0,1)\times(\eta+((p-1)b-1)h,2\eta+(pb-1)h), \quad p=2,\dots,m,\\
D_p^{L} &:= (0,1)\times(0,1)\times((p-1)bh,\eta+pbh), \quad p=1,\dots,m-1,\\
\partial^{L} D_p^{R} &:= (0,1)\times(0,1)\times\{\eta+((p-1)b-1)h\}, \quad p=2,\dots,m,\\
\partial^{R} D_p^{L} &:= (0,1)\times(0,1)\times\{\eta+pbh\}, \quad  p=1,\dots,m-1, \\
X_q^{M} &:= \{(i_1h,i_2h,i_3h):1\le i_1\le n,1\le i_2\le n,(q-1)b+1\le i_3 \le 2\gamma+qb-1\}, \quad q=1,\dots,m, \\
X_p^{R} &:= \{(i_1h,i_2h,i_3h):1\le i_1\le n,1\le i_2\le n,\gamma+(p-1)b\le i_3\le 2\gamma+pb-2\}, \quad p=2,\dots,m, \\
X_p^{L} &:= \{(i_1h,i_2h,i_3h):1\le i_1\le n,1\le i_2\le n,(p-1)b+1\le i_3\le \gamma+pb-1\}, \quad p=1,\dots,m-1. 
\end{align*}

\paragraph{Auxiliary problems.}
For each $q=1,\dots,m$, $H_q^{M} \pmb v=\pmb g$ is defined as the
discretization on $X_q^{M}$ of
\begin{align*}
  &\begin{dcases}
     \left(\partial_1^2+\partial_2^2+(s_q^{M}(x_3)\partial_3)^2+\dfrac{\omega^2}{c^2(x)}\right)v(x)=g(x),&\quad \forall x\in D_q^{M},\\
     v(x)=0,&\quad \forall x\in \partial D_q^{M}, 
   \end{dcases}
\end{align*}
For $p=2,\dots,m$, $H_p^{R} \pmb v=\pmb g$ is defined as the
discretization on $X_p^{R}$ of
\begin{align*}
  &\begin{dcases}
     \left(\partial_1^2+\partial_2^2+(s_p^{R}(x_3)\partial_3)^2+\dfrac{\omega^2}{c^2(x)}\right)v(x)=0,\quad &\forall x\in D_p^{R},\\
     v(x)=w(x_1,x_2),\quad &\forall x\in \partial^{L} D_p^{R}, \\
     v(x)=0,\quad &\forall x\in \partial D_p^{R} \setminus \partial^{L} D_p^{R}, 
   \end{dcases}
\end{align*}
where $\pmb
g:=(-1/h^2)[\pmb w^T,0,\dots,0]^T$ and $\pmb w := [w_{1,1},\dots,w_{n,1},\dots,w_{n,n}]$ is the discrete boundary value. Finally,
for $p=1,\dots,m-1$, $H_p^{L} \pmb v=\pmb g$ is the discretization on
$X_p^{L}$ of
\begin{align*}
  &\begin{dcases}
     \left(\partial_1^2+\partial_2^2+(s_p^{L}(x_3)\partial_3)^2+\dfrac{\omega^2}{c^2(x)}\right)v(x)=0,\quad &\forall x\in D_p^{L},\\
     v(x)=w(x_1,x_2),\quad &\forall x\in \partial^{R} D_p^{L}, \\
     v(x)=0,\quad &\forall x\in \partial D_p^{L} \setminus \partial^{R} D_p^{L},
   \end{dcases}
\end{align*}
where $\pmb g:=(-1/h^2)[0,\dots,0,\pmb w^T]^T$ and $\pmb w := [w_{1,1},\dots,w_{n,1},\dots,w_{n,n}]$. 

\paragraph{Auxiliary Green's operators.}
For $q=1,\dots,m$, $\td{G}_q^{M}:\pmb y\mapsto \pmb z$ is defined
using the following operations:
\begin{enumerate}
\item
  Introduce a vector $\pmb g$ defined on $X_q^{M}$ by setting $\pmb y$ to $X_q$ and zero everywhere else. 
\item
  Solve $H_q^{M} \pmb v=\pmb g$ on $X_q^{M}$. 
\item
  Set $\pmb z$ as the restriction of $\pmb v$ on $X_q$. 
\end{enumerate}
For $p=2,\dots,m$, $\td{G}_p^{R}:\pmb w\mapsto \pmb z$ is given by:
\begin{enumerate}
\item
  Set $\pmb g=(-1/h^2)[\pmb w^T,0,\dots,0]^T$. 
\item
  Solve $H_p^{R} \pmb v =\pmb g$ on $X_p^{R}$. 
\item
  Set $\pmb z$ as the restriction of $\pmb v$ on $X_p$. 
\end{enumerate}
Finally, for $p=1,\dots,m-1$, the operators $\td{G}_p^{L}:\pmb
w\mapsto \pmb z$ is introduced to be:
\begin{enumerate}
\item
  Set $\pmb g=(-1/h^2)[0,\dots,0,\pmb w^T]^T$. 
\item
  Solve $H_p^{L} \pmb v =\pmb g$ on $X_p^{L}$. 
\item
  Set $\pmb z$ as the restriction of $\pmb v$ on $X_p$. 
\end{enumerate}

\paragraph{Putting together.}
In the 3D case, $\pmb g^{L}$ and $\pmb g^{R}$ for the column-major ordering array\\
$\pmb
g=[g_{1,1,1},\dots,g_{s_1,1,1},\dots,g_{s_1,s_2,1},\dots,g_{s_1,s_2,s_3}]^T$
induced from some 3D grid with size\\
$s_1\times s_2 \times s_3$ are
given by
\begin{align*}
  \pmb g^{L}:=[g_{1,1,1},\dots,g_{s_1,1,1},\dots,g_{s_1,s_2,1}]^T, \quad \pmb g^{R}:=[g_{1,1,s_3},\dots,g_{s_1,1,s_3},\dots,g_{s_1,s_2,s_3}]^T. 
\end{align*}

The subproblems $H_q^{M} \pmb v=\pmb g$, $H_p^{R} \pmb v=\pmb g$ and
$H_p^{L} \pmb v=\pmb g$ are quasi-2D. To solve them, we group the
elements along dimension 3 first, and then apply the nested dissection
method\cite{george1973nested,duff1983multifrontal} to them, as in \cite{sweeppml}. This gives the setup process of the 3D
preconditioner in Algorithm \ref{alg:3dsetup} and the application
process in Algorithm \ref{alg:3dapp}.
\begin{algorithm}[h!]
  \caption{Construction of the 3D additive sweeping preconditioner of
    the system \eqref{eqn:3D}. Complexity
    $=O(n^4(b+\gamma)^3/b)=O(N^{4/3}(b+\gamma)^3/b)$.}
  \label{alg:3dsetup}
  \begin{algorithmic}
    \FOR {$q=1,\dots,m$}
    \STATE Construct the nested dissection factorization of $H_q^{M}$, which defines $\td{G}_q^{M}$.
    \ENDFOR
    \FOR {$p=2,\dots,m$}
    \STATE Construct the the nested dissection factorization of $H_p^{R}$, which defines $\td{G}_p^{R}$.
    \ENDFOR
    \FOR {$p=1,\dots,m-1$}
    \STATE Construct the the nested dissection factorization of $H_p^{L}$, which defines $\td{G}_p^{L}$.
    \ENDFOR
  \end{algorithmic}
\end{algorithm}

\begin{algorithm}[h!]
  \caption{Computation of $\td{\pmb u}\approx G \pmb f$ using the
    preconditioner from Algorithm \ref{alg:3dsetup}. Complexity
    $=O(n^3\log n (b+\gamma)^2/b)=O(N\log N(b+\gamma)^2/b)$.}
  \label{alg:3dapp}
  \begin{algorithmic}
    \FOR {$q=1,\dots,m$}
    \STATE
    $\td{\pmb u}_{q,q}=\td{G}_q^{M} \pmb f_q$
    \ENDFOR
    \FOR {$p=2,\dots,m$}
    \STATE
    $\td{\pmb u}_{p,1:p-1}=\td{G}_{p}^{R}\td{\pmb u}_{p-1,1:p-1}^{R}$\\
    $\td{\pmb u}_{p,1:p}^{R}=\td{\pmb u}_{p,1:p-1}^{R}+\td{\pmb u}_{p,p}^{R}$
    \ENDFOR
    \FOR {$p=m-1,\dots,1$}
    \STATE
    $\td{\pmb u}_{p,p+1:m}=\td{G}_p^{L}\td{\pmb u}_{p+1,p+1:m}^{L}$\\
    $\td{\pmb u}_{p,p:m}^{L}=\td{\pmb u}_{p,p}^{L}+\td{\pmb u}_{p,p+1:m}^{L}$
    \ENDFOR
    \STATE
    $\td{\pmb u}_1=\td{\pmb u}_{1,1}+\td{\pmb u}_{1,2:m}$\\
    \FOR {$p=2,\dots,m-1$}
    \STATE
    $\td{\pmb u}_p=\td{\pmb u}_{p,1:p-1}+\td{\pmb u}_{p,p}+\td{\pmb u}_{p,p+1:m}$
    \ENDFOR
    \STATE
    $\td{\pmb u}_m=\td{\pmb u}_{m,1:m-1}+\td{\pmb u}_{m,m}$ 
  \end{algorithmic}
\end{algorithm}

For the algorithm analysis, we notice that each quasi-2D subproblem
has $O(\gamma+b)$ layers along the third dimension. Therefore, the
setup cost for each subproblem is $O((\gamma+b)^3n^3)$ and the
application cost is $O((\gamma+b)^2n^2\log n)$. Taking the total
number of subproblems into account, the total setup cost for the 3D
preconditioner is $O(n^4(b+\gamma)^3/b)$ and the total application
cost is $O(n^3\log n (b+\gamma)^2/b)$.

\subsection{Numerical results}
\label{sec:3Dnumerical}

Here we present the numerical results in 3D. All the settings and
notations are kept the same with Section \ref{sec:2Dnumerical}
unless otherwise stated. The PMLs are put on all sides of the boundary
and the symmetric version of the equation is adopted to save memory
cost. The PML width is $\eta= 9h$ for the boundary and is
$\eta_\text{aux}=5h$ for the interior auxiliary ones. The number of
layers in each subdomain is $b=4$ for the interior ones,
$b+\gamma-1=12$ for the leftmost one and $b+\gamma-2=11$ for the
rightmost one.

\begin{figure}[h!]
  \centering
  \subfigure[]
  {\includegraphics[width=0.32\textwidth]{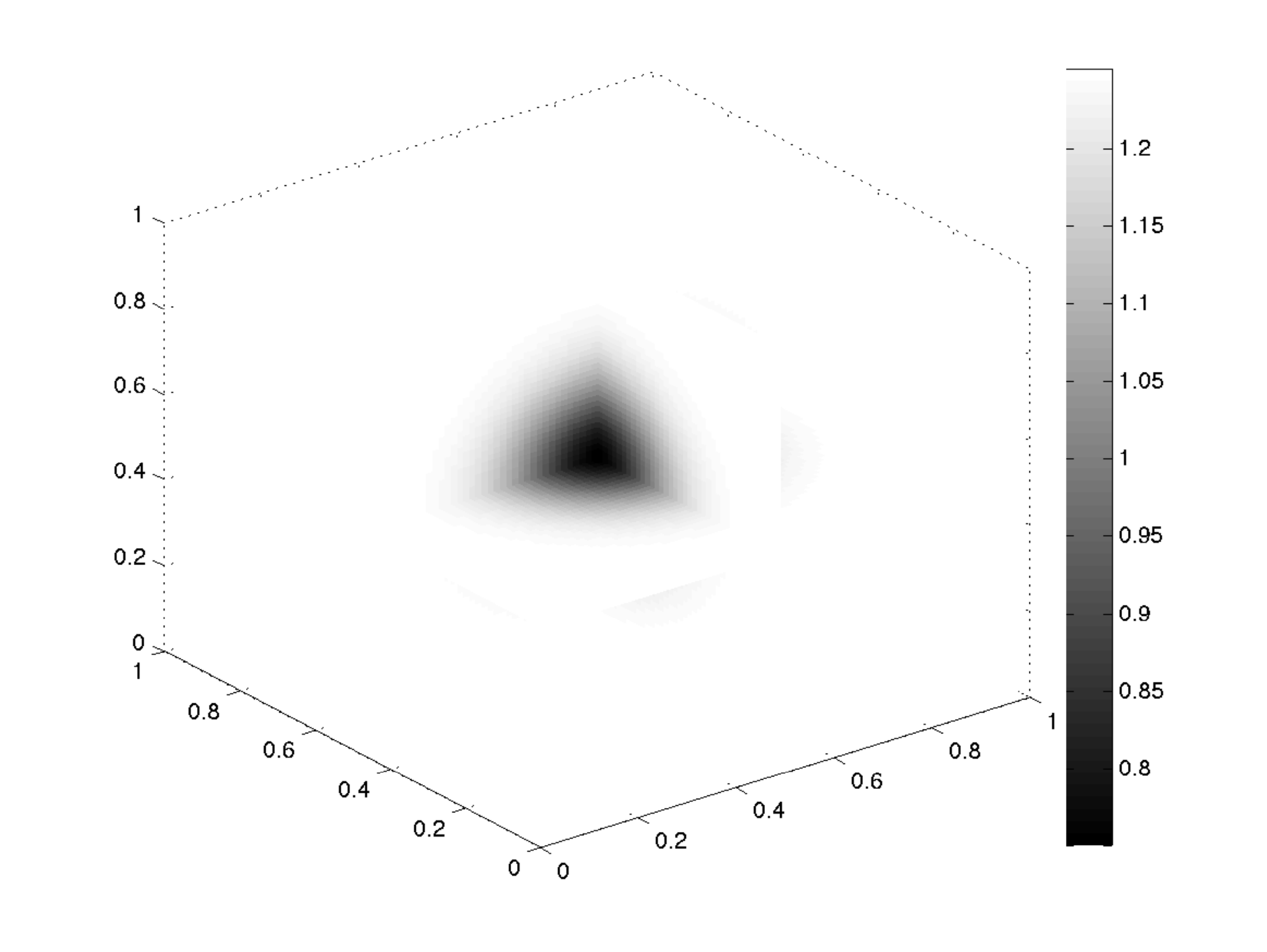}}
  \subfigure[]
  {\includegraphics[width=0.32\textwidth]{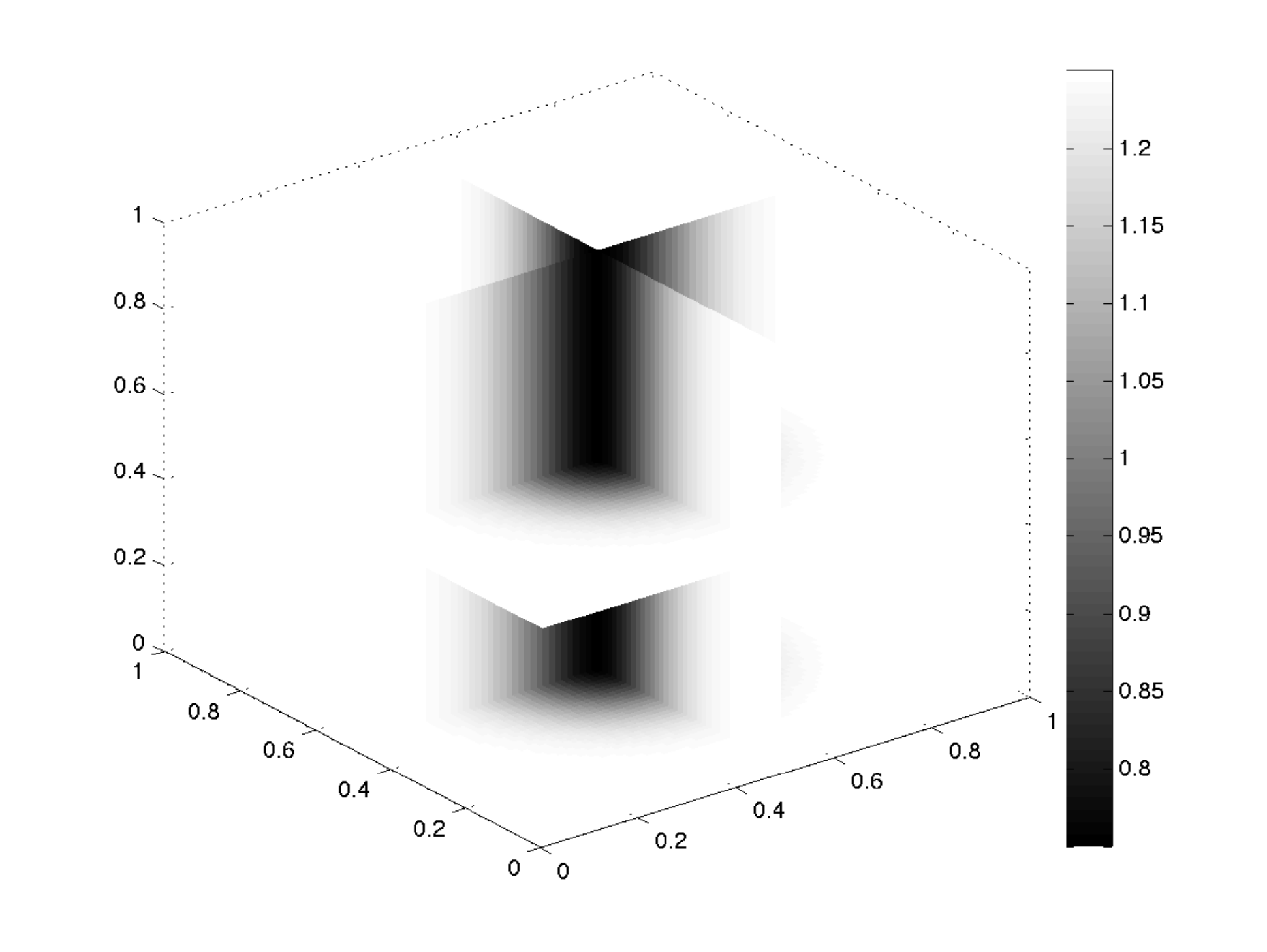}}
  \subfigure[]
  {\includegraphics[width=0.32\textwidth]{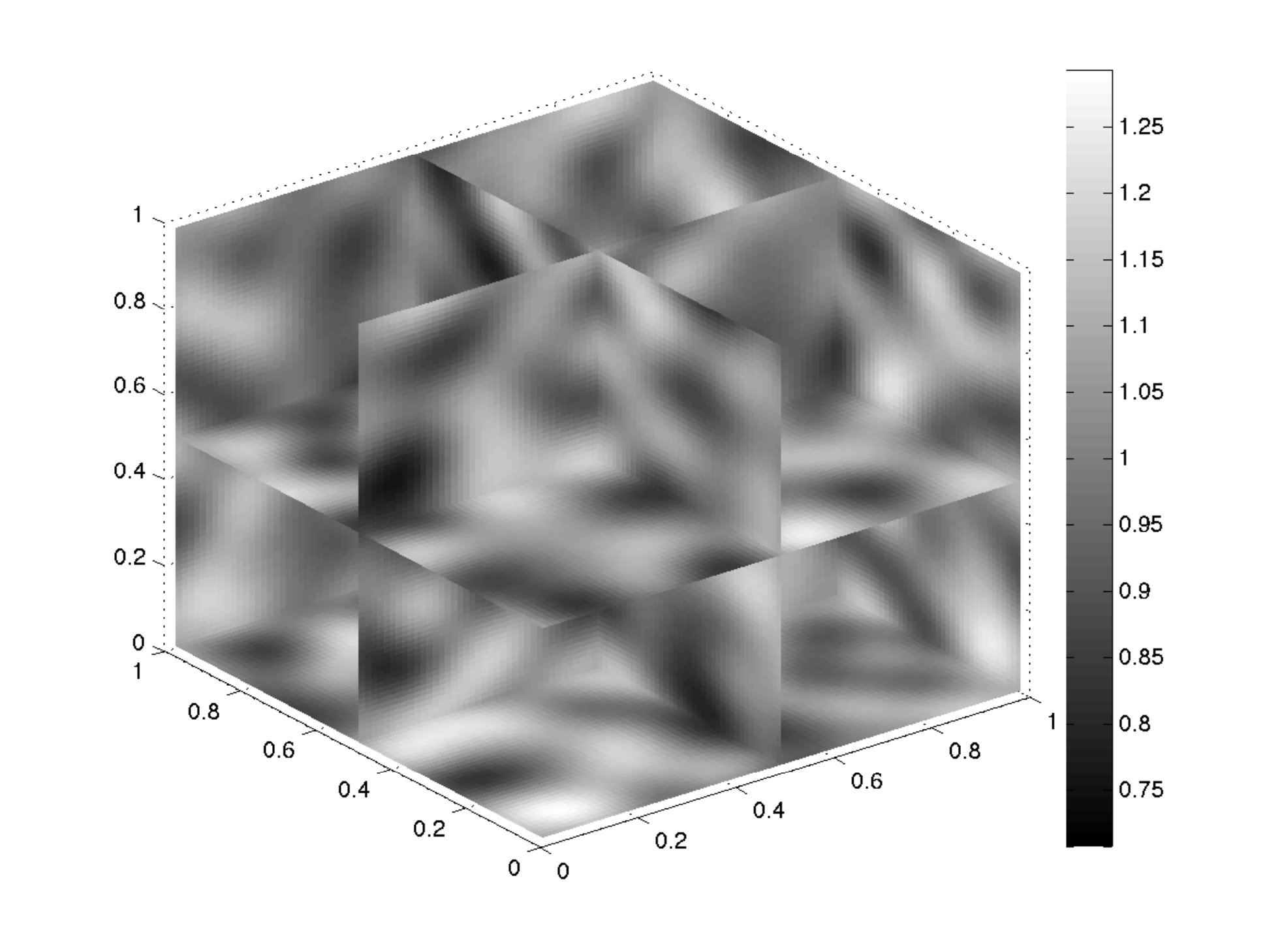}}
  \caption{The three velocity fields tested in 3D.}
  \label{fig:3D}
\end{figure}

The velocity fields tested are (see Figure \ref{fig:3D}):  
\begin{enumerate}[(a)]
\item
  A converging lens with a Gaussian profile at the center of the domain.
\item
  A vertical waveguide with a Gaussian cross-section.
\item
  A random velocity field.
\end{enumerate}

The forces tested for each velocity field are: 
\begin{enumerate}[(a)]
\item
  A Gaussian point source centered at $(1/2,1/2,1/4)$.
\item
  A Gaussian wave packet with wavelength comparable to the typical
  wavelength of the domain. The packet centers at $(1/2,1/4,1/4)$ and
  points to the direction $(0,1/\sqrt{2},1/\sqrt{2})$.
\end{enumerate}

\begin{table}[h!]
\centering
\begin{overpic}
[width=0.45\textwidth]{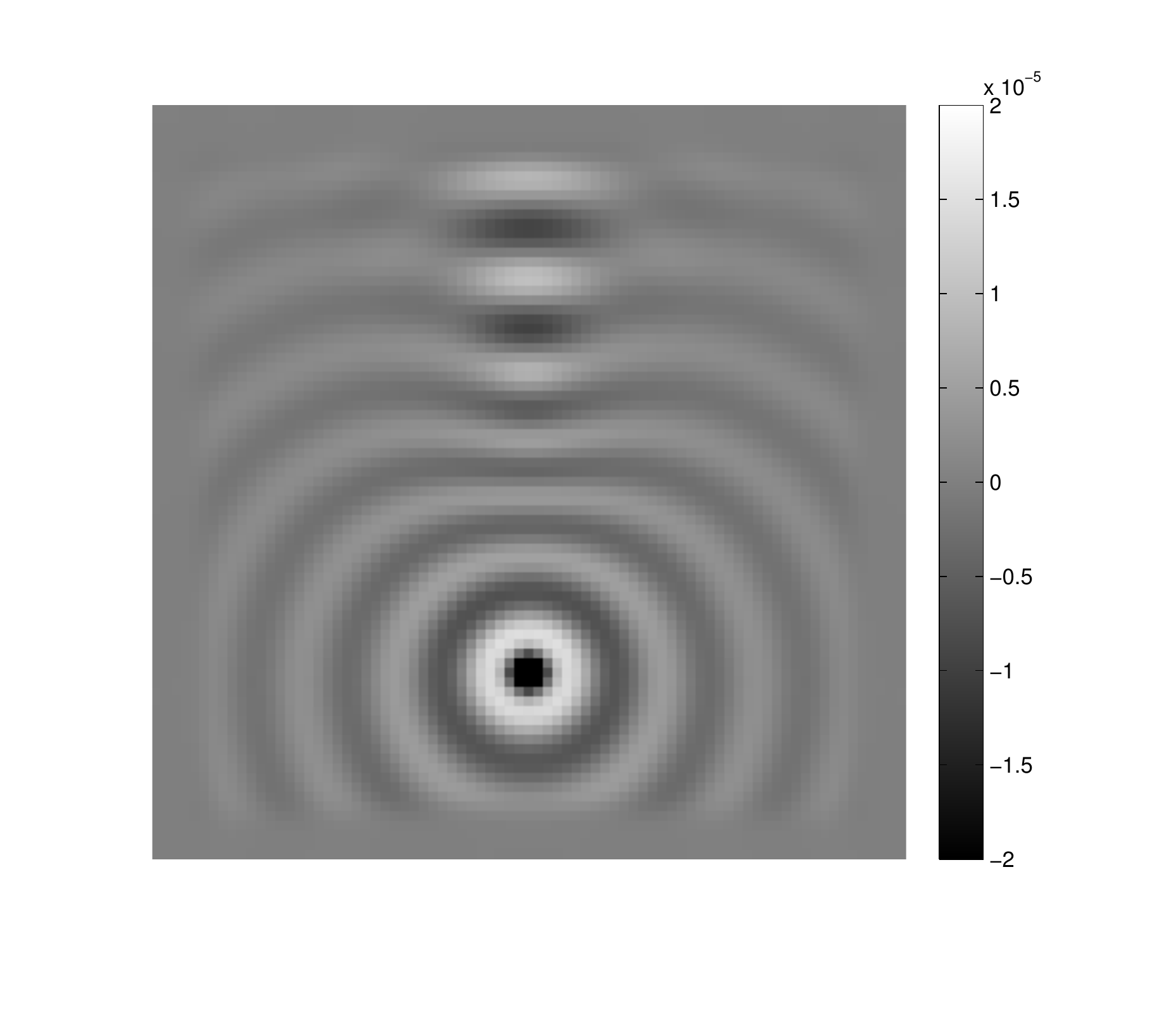}
\put(40,6){force (a)}
\end{overpic}
\begin{overpic}
[width=0.45\textwidth]{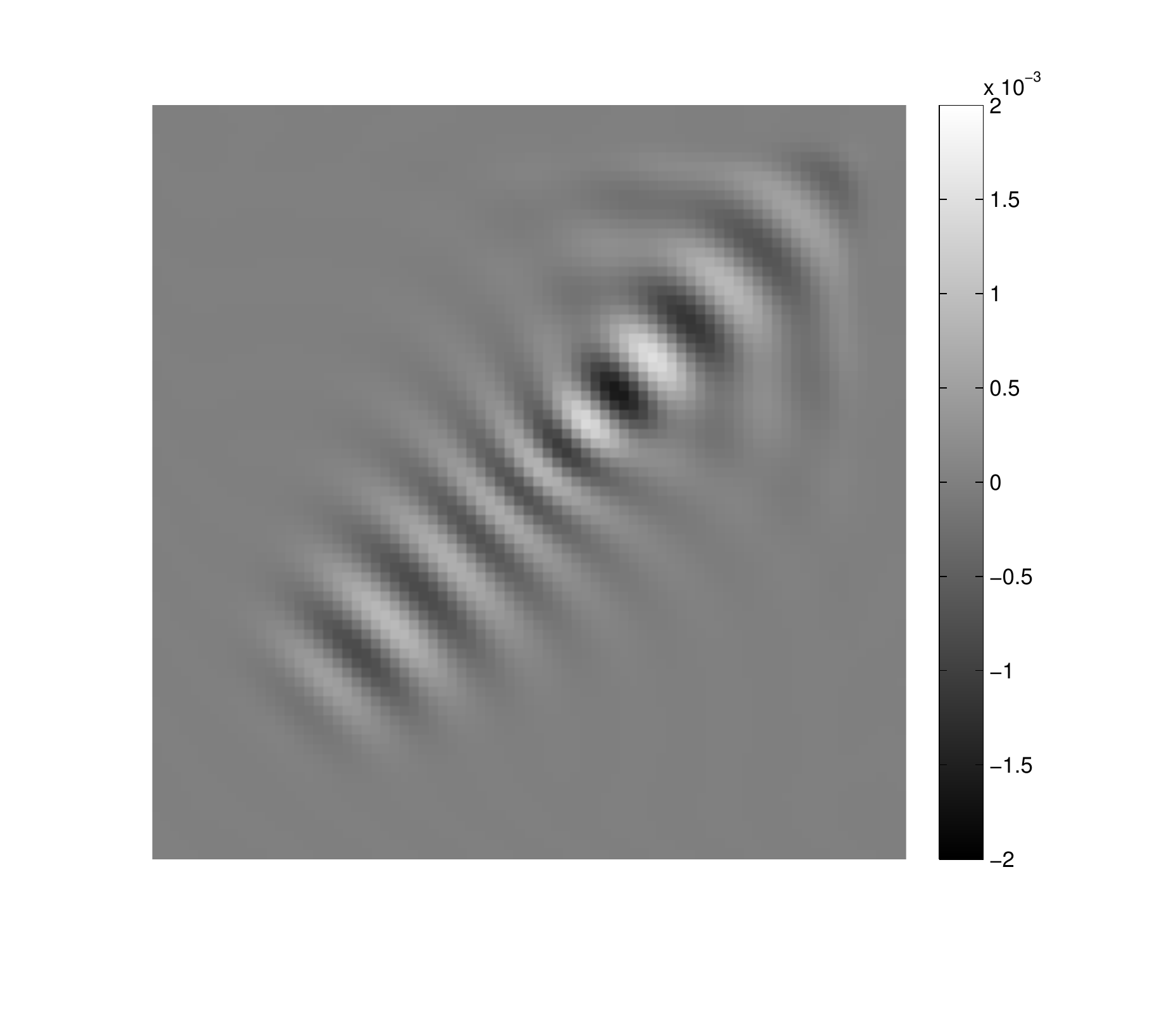}
\put(40,6){force (b)}
\end{overpic}
\begin{tabular}{ccc|cc|cc}
\hline 
\multicolumn{3}{c|}{velocity field (a)} & \multicolumn{2}{c|}{force (a)} & \multicolumn{2}{c}{force (b)} \\ 
\hline 
$\omega/(2\pi)$ & $N$ & $T_{\text{setup}}$ & $N_\text{iter}$ & $T_\text{solve}$ & $N_\text{iter}$ & $T_\text{solve}$ \\ 
\hline 
5 & $39^3$ & 2.3304e$+$01  &3& 2.9307e$+$00  &4&  3.7770e$+$00\\
10 & $79^3$ & 3.2935e$+$02 &3&  3.6898e$+$01  &4&  4.6176e$+$01\\ 
20 & $159^2$ & 4.2280e$+$03  &4&  4.3999e$+$02  &4&  4.6941e$+$02\\
\hline 
\end{tabular} 
\caption{Results for velocity field (a) in 3D. Solutions with $\omega/(2\pi)=10$ at $x_1=0.5$ are presented.}
\label{tab:3domegaa}
\end{table}
\begin{table}[h!]
\centering
\begin{overpic}
[width=0.45\textwidth]{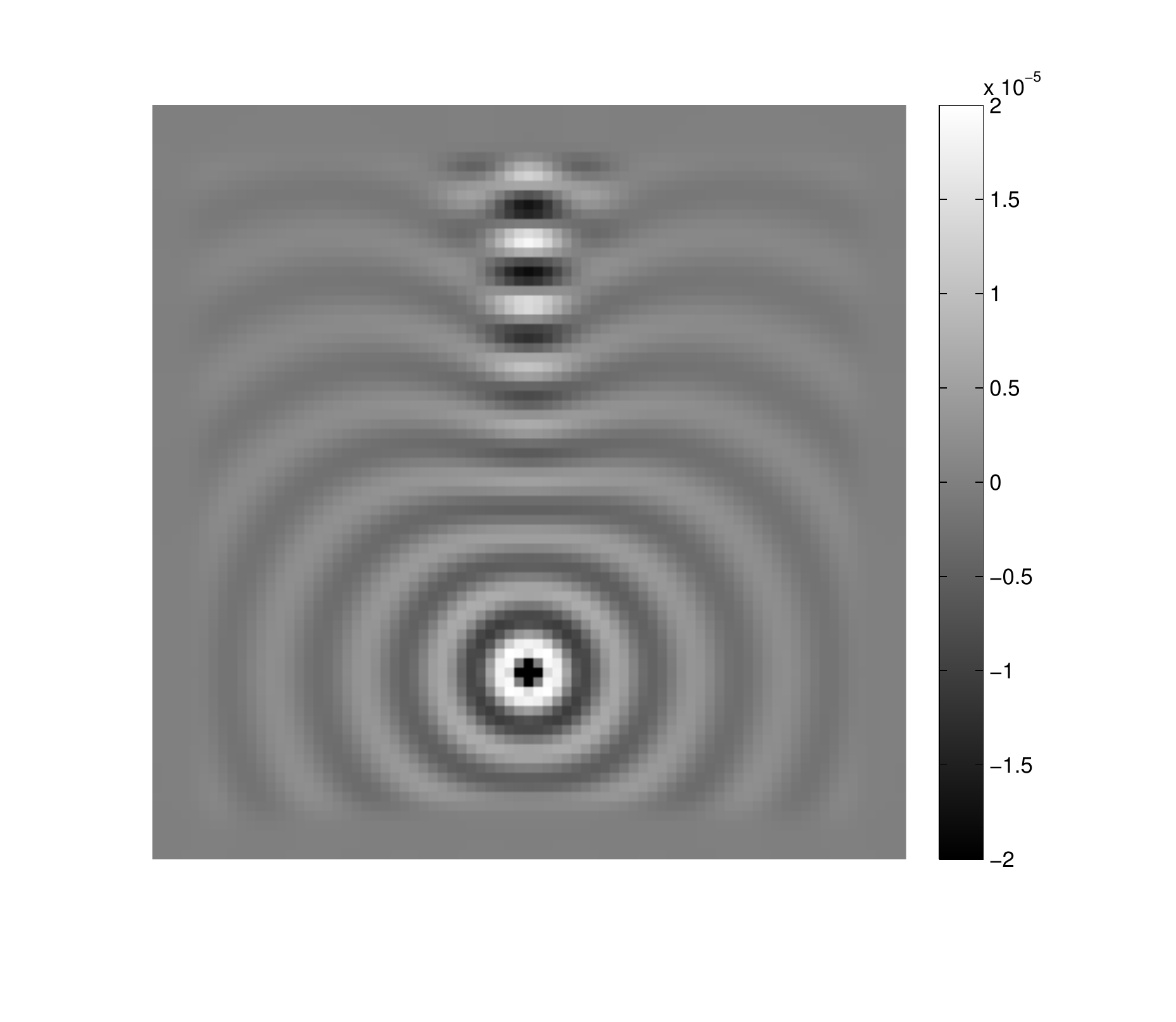}
\put(40,6){force (a)}
\end{overpic}
\begin{overpic}
[width=0.45\textwidth]{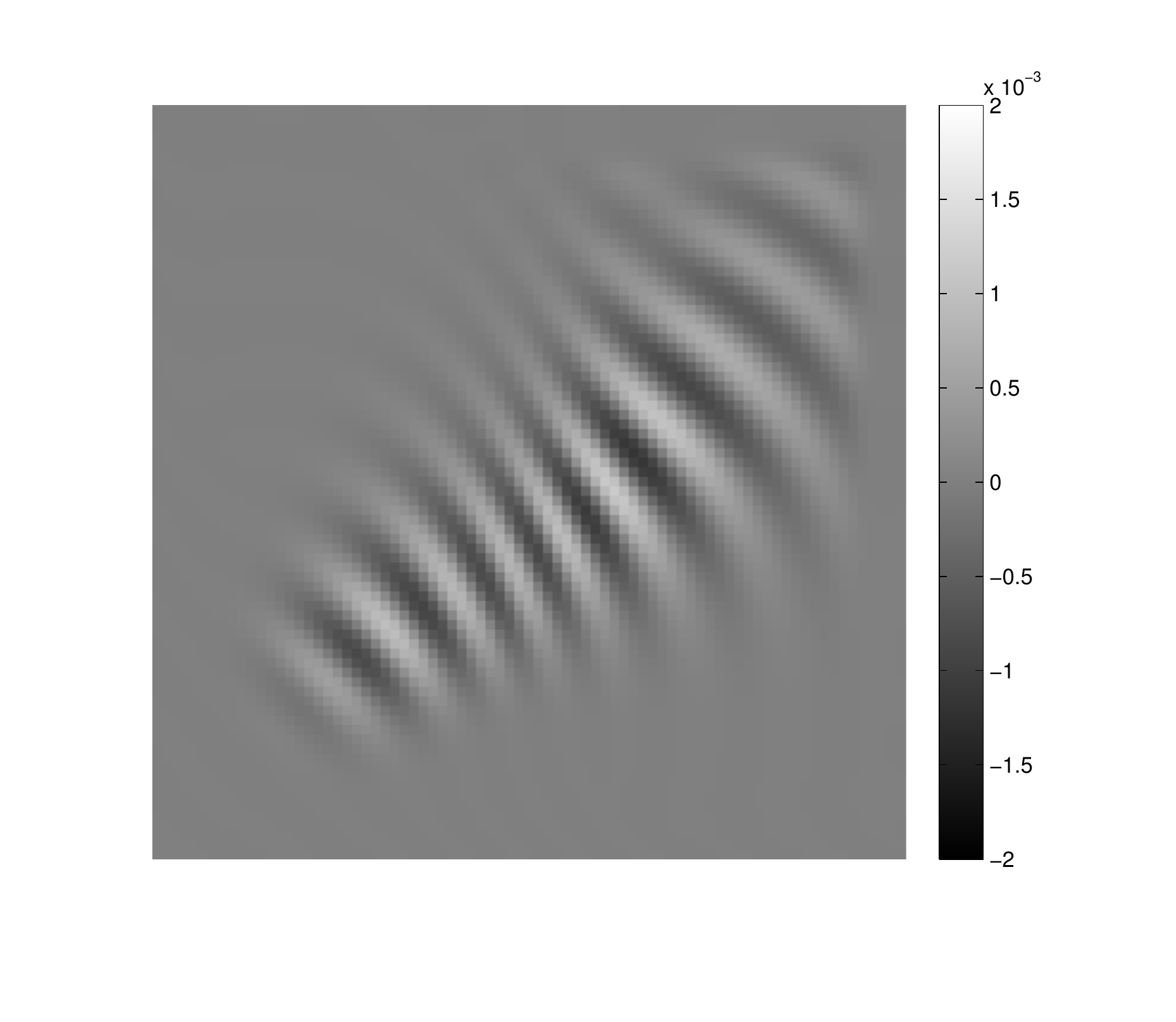}
\put(40,6){force (b)}
\end{overpic}
\begin{tabular}{ccc|cc|cc}
\hline 
\multicolumn{3}{c|}{velocity field (b)} & \multicolumn{2}{c|}{force (a)} & \multicolumn{2}{c}{force (b)} \\ 
\hline 
$\omega/(2\pi)$ & $N$ & $T_{\text{setup}}$ & $N_\text{iter}$ & $T_\text{solve}$ & $N_\text{iter}$ & $T_\text{solve}$ \\ 
\hline 
5 & $39^3$ & 2.1315e$+$01  &3&  2.7740e$+$00  &3&  2.7718e$+$00\\
10 & $79^3$ & 3.4256e$+$02 &4&  4.4286e$+$01  &3&  3.4500e$+$01\\ 
20 & $159^2$ & 4.3167e$+$03  &5&  5.7845e$+$02  &4& 4.6462e$+$02\\
\hline 
\end{tabular} 
\caption{Results for velocity field (b) in 3D. Solutions with $\omega/(2\pi)=10$ at $x_1=0.5$ are presented.}
\label{tab:3domegab}
\end{table}
\begin{table}[h!]
\centering
\begin{overpic}
[width=0.45\textwidth]{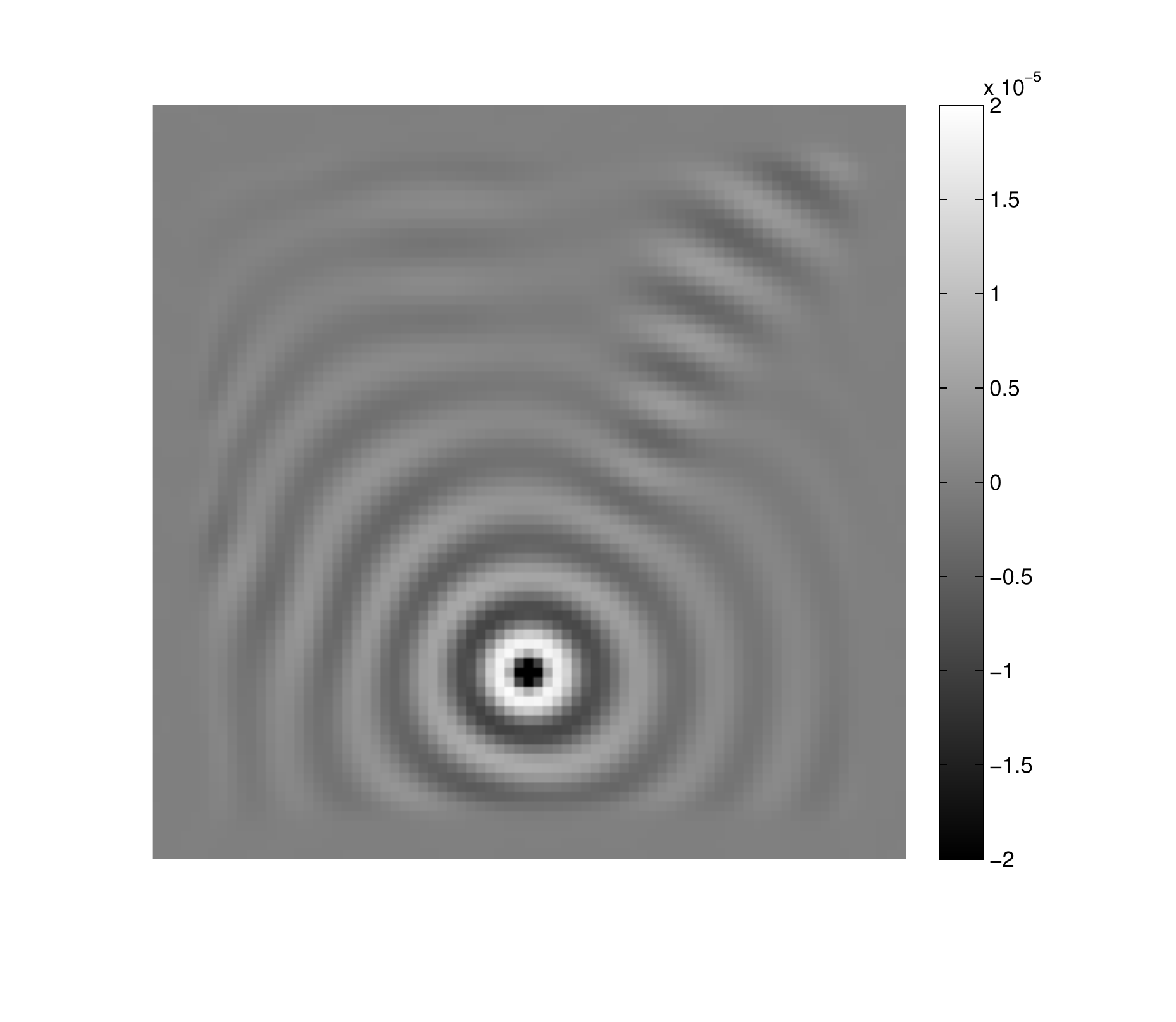}
\put(40,6){force (a)}
\end{overpic}
\begin{overpic}
[width=0.45\textwidth]{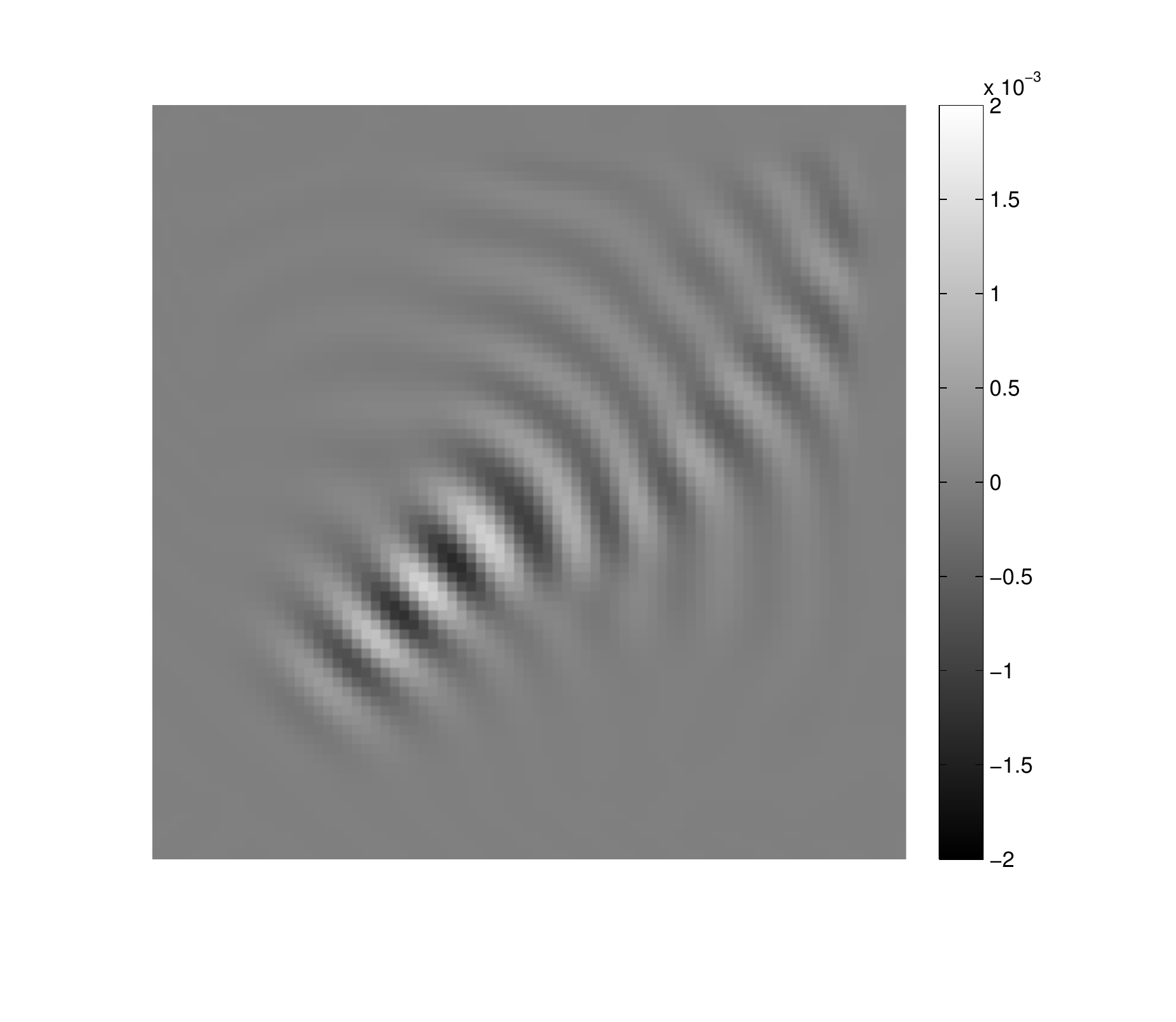}
\put(40,6){force (b)}
\end{overpic}
\begin{tabular}{ccc|cc|cc}
\hline 
\multicolumn{3}{c|}{velocity field (c)} & \multicolumn{2}{c|}{force (a)} & \multicolumn{2}{c}{force (b)} \\ 
\hline 
$\omega/(2\pi)$ & $N$ & $T_{\text{setup}}$ & $N_\text{iter}$ & $T_\text{solve}$ & $N_\text{iter}$ & $T_\text{solve}$ \\ 
\hline 
5 & $39^3$ & 2.1063e$+$01 &4& 3.8074e$+$00  &4&  3.7975e$+$00\\
10 & $79^3$ & 3.4735e$+$02  &4&  4.4550e$+$01  &4&  4.5039e$+$01\\ 
20 & $159^2$ & 4.3391e$+$03  &4&  4.4361e$+$02  &5&  5.8090e$+$02\\
\hline 
\end{tabular} 
\caption{Results for velocity field (c) in 3D. Solutions with $\omega/(2\pi)=10$ at $x_1=0.5$ are presented.}
\label{tab:3domegac}
\end{table}

The results are given in Tables \ref{tab:3domegaa}, \ref{tab:3domegab}
and \ref{tab:3domegac}. From these tests we see that the iteration
number grows mildly as the problem size grows. We also notice that the
setup cost scales even better than $O(N^{4/3})$, mainly because MATLAB
performs dense linear algebra operations in a parallel way, which
gives some extra advantages to the nested dissection algorithm as the
problem size grows.

\section{Conclusion}
\label{sec:Conclusion}
In this paper, we proposed a new additive sweeping preconditioner for
the Helmholtz equation based on the PML. When combined with the
standard GMRES solver, the iteration number grows mildly as the
problem size grows. The novelty of this approach is that the unknowns
are split in an additive way and the boundary values of the
intermediate results are utilized directly. The disadvantage is that,
for each subdomains, three subproblems need to be built up, which is
time consuming compared to \cite{sweeppml} and
\cite{stolk2013domaindecomp}. However, the costly parts of the
algorithm, i.e. the whole setup process and the solve processes of the
subproblems $H_q^{M} \pmb v=\pmb g$, can be done in parallel. The only
parts that must be implemented sequentially are the accumulations of
the left-going and right-going waves, where only the solve processes
of the subproblems $H_p^{L} \pmb v=\pmb g$ and $H_p^{R} \pmb v=\pmb g$
are involved, which are the cheapest parts of the algorithm. Besides,
we think that the whole approximation process is simple and
structurally clear from a physics point of view and the idea might be
easy to be generalized to other equations.

There are also some other directions to make potential
improvements. First, other numerical schemes of the equation and other
approximations of the Sommerfeld radiation condition can be used to
develop more efficient versions of this additive
preconditioner. Second, the parallel version of the nested dissection
algorithm can be combined to solve large scale problems. Last, in the
3D case, the quasi-2D subproblems can be solved recursively by
sweeping along the $x_2$ direction with the same technique, which
reduces the theoretical setup cost to $O(N)$ and the application cost
to $O(N)$. However, compared to \cite{sweeppml}, the coefficient of
the complexity in this new method is larger, so it is not clear
whether or not the recursive approach will be more efficient
practically. Nevertheless, it is of great theoretical interest to look
into it.

\bibliographystyle{abbrv}
\bibliography{references}

\begin{thebibliography}{10}

\bibitem{berenger1994pml}
J.-P. Berenger.
\newblock A perfectly matched layer for the absorption of electromagnetic
  waves.
\newblock {\em J. Comput. Phys.}, 114(2):185--200, 1994.

\bibitem{chen2013sourcetrans}
Z.~Chen and X.~Xiang.
\newblock A source transfer domain decomposition method for {H}elmholtz
  equations in unbounded domain.
\newblock {\em SIAM J. Numer. Anal.}, 51(4):2331--2356, 2013.

\bibitem{chen2013sourcetrans2}
Z.~Chen and X.~Xiang.
\newblock A source transfer domain decomposition method for {H}elmholtz
  equations in unbounded domain {P}art {II}: {E}xtensions.
\newblock {\em Numer. Math. Theory Methods Appl.}, 6(3):538--555, 2013.

\bibitem{chew1994pml}
W.~C. Chew and W.~H. Weedon.
\newblock A {3D} perfectly matched medium from modified {M}axwell's equations
  with stretched coordinates.
\newblock {\em Microw. Opt. Techn. Let.}, 7(13):599--604, 1994.

\bibitem{duff1983multifrontal}
I.~S. Duff and J.~K. Reid.
\newblock The multifrontal solution of indefinite sparse symmetric linear
  equations.
\newblock {\em ACM Trans. Math. Software}, 9(3):302--325, 1983.

\bibitem{sweephmf}
B.~Engquist and L.~Ying.
\newblock Sweeping preconditioner for the {H}elmholtz equation: {H}ierarchical
  matrix representation.
\newblock {\em Comm. Pure Appl. Math.}, 64(5):697--735, 2011.

\bibitem{sweeppml}
B.~Engquist and L.~Ying.
\newblock Sweeping preconditioner for the {H}elmholtz equation: {M}oving
  perfectly matched layers.
\newblock {\em Multiscale Model. Simul.}, 9(2):686--710, 2011.

\bibitem{advances}
Y.~A. Erlangga.
\newblock Advances in iterative methods and preconditioners for the {H}elmholtz
  equation.
\newblock {\em Arch. Comput. Methods Eng.}, 15(1):37--66, 2008.

\bibitem{why}
O.~G. Ernst and M.~J. Gander.
\newblock Why it is difficult to solve {H}elmholtz problems with classical
  iterative methods.
\newblock In {\em Numerical Analysis of Multiscale Problems}, volume~83 of {\em
  Lect. Notes Comput. Sci. Eng.}, pages 325--363. Springer, Heidelberg, 2012.

\bibitem{ailu}
M.~J. Gander and F.~Nataf.
\newblock {AILU} for {H}elmholtz problems: {A} new preconditioner based on the
  analytic parabolic factorization.
\newblock {\em J. Comput. Acoust.}, 9(4):1499--1506, 2001.

\bibitem{george1973nested}
A.~George.
\newblock Nested dissection of a regular finite element mesh.
\newblock {\em SIAM J. Numer. Anal.}, 10(2):345--363, 1973.

\bibitem{johnson2008pmlnotes}
S.~G. Johnson.
\newblock Notes on {P}erfectly {M}atched {L}ayers ({PML}s).
\newblock {\em Lecture notes, Massachusetts Institute of Technology,
  Massachusetts}, 2008.

\bibitem{Liu2015}
F.~{Liu} and L.~{Ying}.
\newblock {Recursive Sweeping Preconditioner for the 3D Helmholtz Equation}.
\newblock {\em ArXiv e-prints}, Feb. 2015.

\bibitem{parallelsweep}
J.~Poulson, B.~Engquist, S.~Li, and L.~Ying.
\newblock A parallel sweeping preconditioner for heterogeneous 3{D} {H}elmholtz
  equations.
\newblock {\em SIAM J. Sci. Comput.}, 35(3):C194--C212, 2013.

\bibitem{stolk2013domaindecomp}
C.~C. Stolk.
\newblock A rapidly converging domain decomposition method for the {H}elmholtz
  equation.
\newblock {\em J. Comput. Phys.}, 241(0):240 -- 252, 2013.

\bibitem{sweepemfem}
P.~Tsuji, B.~Engquist, and L.~Ying.
\newblock A sweeping preconditioner for time-harmonic {M}axwell's equations
  with finite elements.
\newblock {\em J. Comput. Phys.}, 231(9):3770--3783, 2012.

\bibitem{sweepspectral}
P.~Tsuji, J.~Poulson, B.~Engquist, and L.~Ying.
\newblock Sweeping preconditioners for elastic wave propagation with spectral
  element methods.
\newblock {\em ESAIM Math. Model. Numer. Anal.}, 48(2):433--447, 2014.

\bibitem{sweepem}
P.~Tsuji and L.~Ying.
\newblock A sweeping preconditioner for {Y}ee's finite difference approximation
  of time-harmonic {M}axwell's equations.
\newblock {\em Front. Math. China}, 7(2):347--363, 2012.

\bibitem{vion2014doublesweep}
A.~Vion and C.~Geuzaine.
\newblock Double sweep preconditioner for optimized {S}chwarz methods applied
  to the {H}elmholtz problem.
\newblock {\em J. Comput. Phys.}, 266(0):171 -- 190, 2014.

\bibitem{demanet}
L.~{Zepeda-N{\'u}{\~n}ez} and L.~{Demanet}.
\newblock The method of polarized traces for the 2{D} {H}elmholtz equation.
\newblock {\em ArXiv e-prints}, Oct. 2014.

\end{thebibliography}
\end{document}